\documentclass[twoside,12pt]{article}

\usepackage{amssymb, amsmath, latexsym, mathrsfs, verbatim, calc, hyperref}

\numberwithin{equation}{section} 
\newtheorem{thm}{Theorem}[section]
\newtheorem{lem}[thm]{Lemma}
\newtheorem{cor}[thm]{Corollary}
\newtheorem{prop}[thm]{Proposition}
\newtheorem{rem}[thm]{Remark}

\newtheorem{defn}[thm]{Definition}

\def\ba{\begin{array}}
\def\ea{\end{array}}
\def\beq{\begin{equation}}
\def\bes{\begin{equation*}}
\def\ees{\end{equation*}}
\def\bea{\begin{eqnarray}}
\def\eea{\end{eqnarray}}
\def\beas{\begin{eqnarray*}}
\def\eeas{\end{eqnarray*}}

\def\dis{\displaystyle}

\def\no{\noindent}

\def\norm{|\nts|}

\def\lastline{\par \vspace{-7.3ex} \no}

\def\nts{\negthinspace}
\def\ss{\smallskip}
\def\ms{\medskip}
\def\bs{\bigskip}
\def\q{\quad}
\def\qq{\qquad}

\def\ua{\mathop{\uparrow}}

\def\Ra {\mathop{\Rightarrow }}

\def\La {\mathop{\Leftarrow }}

\def\={=\nts \nts=\nts \nts=\nts \nts=}

\def\wt{\widetilde}

\def\({\textnormal{(}}
\def\){\textnormal{)}}
\def\[{[\neg[}
\def\]{]\neg]}

\def\cd{\cdot}

\def\fa{\,\forall \,}
\def\pa{\partial}
\def\es{\emptyset}

\def\d{\delta}
\def\e{\varepsilon}
\def\z{\zeta}

\def\l{\lambda}

\def\si{\sigma}

\def\t{\tau}
\def\f{\varphi}

\def\o{\omega}

\def\f{\phi}
\def\vf{\varphi}

\def\D{\Delta}

\def\L{\Lambda}
\def\O{\Omega}

\def\bF{{\bf F}}

\def\cD{{\cal D}}
\def\cE{{\cal E}}
\def\cF{{\cal F}}
\def\cG{{\cal G}}
\def\cH{{\cal H}}

\def\cM{{\cal M}}

\def\hB{\mathbb{B}}
\def\hC{\mathbb{C}}
\def\hD{\mathbb{D}}
\def\hE{\mathbb{E}}

\def\hN{\mathbb{N}}

\def\hR{\mathbb{R}}

\def\sB{\mathscr{B}}

\def\sE{\mathscr{E}}

\def\sN{\mathscr{N}}

\def\sP{\mathscr{P}}

\def\esssup{\mathop{\rm esssup}}

\def\dtp{{\hbox{$dt \times dP$-a.s.}}}
\def\pas{{\hbox{$P$-a.s.}}}

\def\sa{\searrow}

\def\no{\noindent}

\def\ss{\smallskip}
\def\ms{\medskip}
\def\bs{\bigskip}
\def\q{\quad}
\def\qq{\qquad}
\def\hb{\hbox}
\def\pa{\partial}
\def\cd{\cdot}

\def\lan{\langle}
\def\ran{\rangle}

\def\bF{{\bf F}}

\def\neg{\negthinspace}
\def\dneg{\neg \neg}
\def\tneg{\neg \neg \neg}
\def\qed{\hfill \rule[0cm]{.25cm}{.25cm}\medskip}   
\def\dfnn{\stackrel{\triangle}{=}}

\def\lan{\langle}
\def\ran{\rangle}
\def\b1{{\bf 1}}

\newenvironment{itm}{\vspace{-1ex}\begin{itemize}}{\end{itemize}}
\def\bi{\begin{itm}}
\def\ei{\end{itm}}

\def\equ_ind{\arabic{section}.\arabic{equation}}
\def\sec_ind{\arabic{section}}

\begin{document}

\title{\bf On Quadratic $g$-Evaluations/Expectations and Related Analysis}

\author{
Jin Ma,\thanks{ \noindent Department of Mathematics, Purdue
University, West Lafayette, IN 47907-1395; Department of
Mathematics, University of Southern California, 3620 S. Vermont
Ave., KAP 108, Los Angeles, CA 90089. Email: jinma@usc.edu.  This
author is supported in part by NSF grant \#0505427. }\q Song
Yao\thanks{ \noindent Department of Mathematics, University of Michigan,
 Ann Arbor, MI 48109; email: songyao@umich.edu. } }

\date{}

\maketitle

\centerline{\bf Abstract }

\bs

In this paper we extend the notion of $g$-evaluation, in particular
$g$-expectation, of Peng \cite{Peng-97, Pln} to the case where the
generator $g$ is allowed to have a quadratic growth (in the variable
``$z$"). We show that
some important properties of the $g$-expectations, including a
representation theorem between the generator and the corresponding
$g$-expectation,
and consequently the reverse comparison theorem of quadratic BSDEs
as well as the Jensen inequality, remain true in the quadratic case.
Our main results also include a Doob-Meyer type decomposition, the
optional sampling theorem, and the up-crossing inequality.
The results of this paper are important in
the further development of the general quadratic nonlinear
expectations (cf. \cite{HMPY-07}).
 \vfill \vspace{.5cm} \no {\bf Keywords: }\:
Quadratic $g$-evaluations, quadratic $g$-expectations, BMO, reverse
comparison theorem, Jensen's inequality, Doob-Meyer Decomposition,
optional sampling, upcrossing inequality.

\eject

\section{Introduction}
\setcounter{equation}{0}

In this paper we extend the notion of {\it $g$-evaluations},
introduced by  Peng \cite{Pln}, to the case when the generator $g$
is allowed to have quadratic growth in the variable $z$. This will
include the so-called quadratic $g$-expectation as a special case,
as was in the linear growth case initiated in \cite{Peng-97}. The
notion of $g$-expectation, as a nonlinear extension of the
well-known Girsanov transformations and originally motivated by
theory of expected utility, has been found to have direct relations
with a fairly large class of risk measures in finance. When the
nonlinear expectation is allowed to have possible quadratic growth,
it is expected that it will lead to the representation theorem that
characterizes the general convex, but not necessarily ``coherent"
risk measures in terms of a class of quadratic BSDEs. The most
notable example of such risk measure is the entropic risk measure
(see, e.g., Barrieu and El Karoui \cite{BarKar}), which is known to
have a representation as the solution to a quadratic BSDE, but falls
outside the existing theory of the ``filtration-consistent nonlinear
expectations" \cite{CHMP}, which requires that the generator be only
of linear growth. We refer the readers to \cite{Peng-97},
\cite{BCHMP}, \cite{CHMP}, and the expository paper \cite{Pln} for
more detailed account for basic properties of $g$-evaluations and
$g$-expectations, as well as the relationship between the risk
measures and $g$-expectations. A brief review of the basic
properties of $g$-evaluations and $g$-expectations will be given in
\S2 for ready references.

The main purpose of this paper is to introduce the notion of
quadratic $g$-evaluation and $g$-expectation, and prove some of the
important properties that are deemed as essential. In an
accompanying paper \cite{HMPY-07} we shall further extend the notion
of filtration consistent nonlinear expectation to the quadratic
case, and establish the ultimate relations between a convex risk
measure and a BSDE. The main results in this paper include the
Doob-Meyer decomposition theorem, optional sampling theorem,
upcrossing inequality, and Jensen's inequality. We also prove that
the quadratic generator can be represented as the limit of
the difference quotients of the corresponding $g$-evaluation,
extending the result in linear growth case \cite{BCHMP}. With the
help of this result, we can then prove the so-called {\it reversed
comparison theorem}, as in the linear case.

Although most of the results presented in this paper look similar to
those in the linear case, the techniques involved in the proofs are
quite different. We combine the techniques used in the study for
quadratic BSDEs, initiated by Kobylanski \cite{Ko} and the by now
well-known properties of the BMO martingales. Since many of these
results are interesting in their own right, we often present full
details of proofs for future references.

\ms

This paper is organized as follows. In section 2 we give the
preliminaries, and review the existing theory of
$g$-evaluation/expectations and BMO martingales. In section 3 we
define the quadratic $g$-evaluation and discuss its basic
properties. Some fine properties of $g$-evaluations/expectations are
presented in Section 4. These include a representation of quadratic
generator via quadratic $g$-evaluations, a reverse comparison
theorem of quadratic BSDE, and the Jensen's inequality. In section 5
we prove the main results of this paper regarding the quadratic
$g$-martingales: a Doob-Meyer type decomposition, the Optional
Sampling theorem, and the Upcrossing Inequality.

\section{Preliminaries}
\setcounter{equation}{0}

Throughout this paper we consider a filtered, complete probability
space $(\O,\cF, P, \bF)$ on which is defined a $d$-dimensional
Brownian motion $B$. We assume that the filtration $\bF\dfnn
\{\cF_t\}_{t\ge0}$ is generated by the Brownian motion $B$,
augmented by all $P$-null sets in $\cF$, so that it satisfies the
{\it usual hypotheses} (cf. \cite{Pr-90}).  We denote $\sP$ to be
the progressively measurable $\si$-field on $\O\times [0,T]$; and
$\cM_{0,T}$ to be the set of all $\bF$-stopping times $\tau$ such
that $0\leq\tau\leq T$, $P$-a.s., where $T>0$ is some fixed
time horizon. 

In what follows we fix a finite time horizon $T>0$, and denote $\hE$
to be a generic Euclidean space, whose inner product and norm will
be denoted by $\langle\cdot,\cdot\rangle$ and $|\cdot|$,
respectively; and denote $\hB$ to be a generic Banach space with
norm $\|\cd\|$. Moreover,
the following spaces of functions will be frequently used in the
sequel. Let $\cG$  be a generic sub-$\si$-field of $\cF$, we denote

\vspace{-5pt}
\begin{itemize}
\item for $0\le p\le\infty$, $L^p(\cG;\hE)$ to be all $\hE$-valued,
$\cG$-measurable random variables $\xi$, with $E(|\xi|^p)<\infty$.
In particular, if $p=0$,
then $L^0(\cG;\hE)$ denotes the space of all $\hE$-valued,
$\cG$-measurable random variables; and if $p=\infty$, then
$L^\infty(\cG;\hE)$ denotes the space of all $\hE$-valued,
$\cG$-measurable random variables $\xi$ such that $\|\xi\|_\infty
\dfnn \underset{\o \in \O}{\esssup}|\xi(\o)|<\infty$;

\item $0 \le  p\le\infty$, $L^p_\bF([0,T];\hB)$ to be all
$\hB$-valued, $\bF$-adapted processes $\psi$, such that
$E\int_0^T\|\psi_t\|^pdt<\infty$. In particular, $p=0$ stands for
all $\hB$-valued, $\bF$-adapted processes; and $p=\infty$ denotes
all processes $X\in L^0_\bF([0,T];\hB)$
such that $\|X\|_\infty \dfnn \underset{t,\o}
{\esssup} |X(t,\o)|<\infty$;

\item $\hD^\infty_\bF([0,T];\hB)=\{X\in L^\infty_\bF([0,T];\hB): \hb{ $X$ has c\`adl\`ag paths}\}$;

\item $\hC^\infty_\bF([0,T];\hB)=\{X\in \hD^\infty_\bF([0,T];\hB): \hb{ $X$ has continuous
paths}\}$;

\item $\cH^2_\bF([0,T];\hB)=\{X \in L^2_\bF([0,T];\hB): \hb{$X$ is predictably
measurable}\}$.
\end{itemize}

Finally, if $d=1$, we shall drop $\hE=\hR$ from the notation (e.g.,
$L^p_\bF([0,T])=L^p_\bF([0,T];\hR)$,
$L^\infty({\cF_T})=L^\infty({\cF_T};\hR)$, and so on).

\bs

\no {\bf $g$-Evaluations and $g$-Expectations}

\ms

We first recall the notion of $g$-evaluation introduced in Peng
\cite{Pln}. Given a time duration $[0,T]$, and a ``generator"
$g=g(t,\o,y,z): [0,T] \times \O \times \hR \times \hR^d \mapsto \hR$
satisfying the standard conditions (e.g., it is Lipschitz in all
spatial variables, and is of linear growth, etc.), consider the
following BSDE on $[0,t]$, $t\in[0,T]$:
 \bea
 \label{BSDE}
 Y_s=\xi+\int_s^t g(r,Y_r,Z_r)dr - \int_s^t Z_r dB_r,  \qq  s \in
 [0,t],
 \eea
where $\xi \in L^2(\cF_t)$. Denote the unique solution by
$(Y^{t,\xi},Z^{t,\xi})$.
The {\it $g$-evaluation} is defined as the family of operators
$\big\{\cE^g_{s,t}: L^2(\cF_t) \mapsto L^2(\cF_s) \big\}_{0 \le s
\le t \le T}$ such that for any $t \in [0,T]$,
$ \cE^g_{s,t}[\xi]\dfnn Y^{t,\xi}_s$, $s \in [0, t]$.

In particular, for any $\xi \in L^2(\cF_T)$, its {\it
$g$-expectation} is defined by $\cE^g(\xi)\dfnn Y^{T,\xi}_0$, and
its {\it conditional $g$-expectation} is defined by
$\cE^g[\xi|\cF_t] \dfnn \cE^g_{t,T}[\xi]$, for any $t \in [0,T]$. We
shall denote (\ref{BSDE}) by BSDE$(t, \xi,g)$ in the sequel for
notational convenience.

\ms

\begin{rem} {\rm An important ingredient in the definition of
$g$-evaluation is its ``domain", namely the subset in $L^0(\cF_T)$
on which the operator is defined (in the current case being
naturally taken as $L^2(\cF_T)$). The domain of a
$g$-evaluation/expectation may vary as the conditions on the
coefficients change, due to the restrictions on the well-posedness
of the BSDE (\ref{BSDE}). For example, owing to the nature of
quadratic BSDEs, in the rest of this paper we shall choose
$L^\infty(\cF_T)$ as the domain for quadratic $g$-evaluations. We
refer to our accompanying paper \cite{HMPY-07} for a more detailed
discussion on the issue of domains for general nonlinear
expectations. \qed
}
\end{rem}

By virtue of the uniqueness of the solution $(Y^{t,\xi},Z^{t,\xi})$,
one can show that the $g$-evaluation $\cE^g_{s,t}$ has the following
properties:
 \bi
 \item[\(1\)] ({\it Monotonicity}) For any $\xi, \eta \in L^2(\cF_t)$ with $\xi  \ge  \eta $, \pas,    $ \cE^g_{r,t}[\xi]  \ge \cE^g_{r,t}[\eta]$, $\pas$;

\item[\(2\)] ({\it Time-Consistency})~~$\cE^g_{r,s}\big[\cE^g_{s,t}[\xi]\big]
=\cE^g_{r,t}[\xi]$, $\pas$, $\xi \in L^2(\cF_t)$, $0 \le r \le s \le
t \le T$;

\if{0}
\item[\(3\)] ({\it Constant-Preserving})~~$\cE^g_{s,t}[\xi]=\xi$, $\pas$, $\xi \in
L^(\cF_s)$, if it holds $\dtp$ that
 \bea\label{g0}
 g(t,\o,y,0)=0, \q~~~y \in \hR;
 \eea
\fi

\item[\(3\)] ({\it Constant-Preserving})~~$\cE^g_{s,t}[\xi]=\xi$, $\pas$, $\xi \in L^2(\cF_s)$,
if it holds $dt\times dP$-a.s. that
 \bea\label{g0}
 g(t,y,0)=0, \qq y \in \hR;
 \eea

\if{0}
\item[\(4\)] (``{\it Zero-one Law}")~~$ \b1_A\cE^g_{s,t}[\b1_A \xi]=\b1_A
\cE^g_{s,t}[\xi]$, $\pas$, $\fa A \in \cF_s$, $\fa \xi \in \L_t$, if
 \bea \label{ga1}
 \b1_A \xi \in \L_r, \q~~ \fa A \in \cF_r,~ \fa \xi
\in \L_r, ~ \fa r \in [0,T];
 \eea
\fi

\item[\(4\)] (``{\it Zero-one Law}")~~For any $\xi \in
L^2(\cF_t)$ and any $A \in \cF_s$, $s \in [0,t]$, it holds that
 $$ \b1_A \cE^g_{s,t}[\xi]=\b1_A\cE^g_{s,t}[\b1_A \xi], \q\pas
 $$
Moreover, if $g(t,0,0)=0$, $dt\times dP$-a.s., then $\b1_A
\cE^g_{s,t}[\xi]=\cE^g_{s,t}[\b1_A \xi]$, $\pas$;

\if{0}
\item[\(5\)] (Translation Invariance)~~$\cE^g_{s,t}[\xi+\eta]=\cE^g_{s,t}[\xi]+\eta$,
$\pas$, $\xi \in L^2(\cF_t)$, $\eta \in L^2(\cF_s)$, if $g$ is
independent of $y$ and $\L$ is closed under addition.
\fi

\item[\(5\)] ({\it Translation Invariance})~~Assume that $g$ is independent of $y$, then for
any $\xi \in L^2(\cF_t)$ and $\eta \in L^2(\cF_s)$,
it holds that $ \cE^g_{s, t}[\xi+\eta]=\cE^g_{s,t}[\xi]+\eta$,
$\pas$

\ei

Clearly, if $g$ satisfies (\ref{g0}), then one can deduce from (2)
and (3) above
that
 \bea\label{coin}
\cE^g[\xi|\cF_s]=\cE^g_{s,T}[\xi]=\cE^g_{s,t}\big[\cE^g_{t,T}[\xi]\big]=\cE^g_{s,t}[\xi],\q
\pas, ~~~ \xi \in L^2(\cF_t), ~~~ 0 \le s \le t \le T;
 \eea
and the conditional $g$-expectation $\cE^g\{\cdot|\cF_t\}$ possesses
the following properties that more or less justify its name
(assuming (\ref{g0}) for (2a) and (3a) below):
 \bi
  \item[\(1a\)] ({\it Monotonicity}) For any $\xi, \eta \in L^2(\cF_T)$ with $\xi  \ge  \eta $, \pas,    $\cE^g[\xi|\cF_t]   \ge \cE^g[\eta |\cF_t]  $, $\pas$;

\item[\(2a\)] (Time-Consistency)~~$\cE^g\big[\cE^g[\xi|\cF_t]\big|\cF_s\big]=\cE^g[\xi|\cF_s]$,
$\pas$, $\xi \in L^2(\cF_T)$, $s \in [0,t]$;
\item[\(3a\)] (Constant-Preserving)~~$\cE^g[\xi|\cF_t]=\xi$, $\pas$, $\xi \in
L^2(\cF_t)$;
\item[\(4a\)] (Zero-one Law)~~For any $\xi \in
L^2(\cF_T)$ and $A \in \cF_t$, it holds that
$\b1_A \cE^g[\b1_A \xi |\cF_t]=\b1_A \cE^g[ \xi |
\cF_t]$, $\pas$; Moreover, if $g(t,0,0)=0$, $\dtp$, then $\b1_A
\cE^g[\xi|\cF_t]=\cE^g[\b1_A \xi|\cF_t]$, $\pas$;
\item[\(5a\)] (Translation Invariance)~~Assume that $g$ is independent of $y$, then
for any $\xi \in L^2(\cF_T)$ and $\eta \in L^2(\cF_t)$ it holds that
$ \cE^g[\xi+\eta|\cF_t]=\cE^g[\xi|\cF_t]+\eta$, $\pas$
 \ei

\ms

\no {\bf BMO Martingales and BMO Processes}

\ms

An important tool for studying the quadratic BSDEs, whence the
quadratic $g$-expectations, is the so-called ``BMO martingales"
and the related stochastic exponentials (see, e.g., \cite{HIM}). We
refer to the monograph of Kazamaki \cite{Ka} for a complete
exposition of the theory of continuous BMO and exponential
martingales. In what follows we shall list some of the important
facts that are useful in our future discussions for ready
references.

To begin with, we recall that a uniformly integrable martingale $M$
null at zero is called a ``BMO martingale" on $[0,T]$ if for some
$1\le p<\infty$, it holds that
 \bea
 \label{BMOp}
 \|M\|_{BMO_p}\dfnn\sup_{\t\in\cM_{0,T}}\Big\|E\{|M_T-M_{\t-}|^p\big|
 \cF_\t\}^{1/p}\Big\|_\infty<\infty.
 \eea
In such a  case we denote $M\in$BMO$(p)$. It is important to note
that  $M\in$BMO$(p)$ if and only if $M\in$BMO$(1)$, and all the
BMO$(p)$ norms are equivalent (cf. \cite{Ka}). Therefore in what
follows we shall say that a martingale $M$ is BMO without specifying
the index $p$; and we shall use only the BMO$(2)$ norm and denote it
simply by $\|\cd\|_{BMO}$. Note also that for a {\it continuous}
martingale $M$ one has
 \beas
 \|M\|_{BMO}=\|M\|_{BMO_2}=\sup_{\t\in\cM_{0,T}}
\Big\|E\{\lan M\ran_T-\lan M\ran_\t\big|\cF_\t\}^{1/2}\Big\|_\infty.
 \eeas

For a given Brownian motion $B$, we say that a process $Z\in
L^2_\bF([0,T];\hR^d)$ is a BMO process, denoted by $Z \in$ BMO by a
slight abuse of notations, if the stochastic integral $M\dfnn Z\cdot
B=\int Z_tdB_t$ is a BMO martingale.

Next,
for a continuous martingale $M$, the Dol\'eans-Dade stochastic
exponential of $M$, denoted customarily by $\sE(M)$, is defined as
$\sE(M)_t\dfnn \exp\{M_t-\frac12\lan M\ran_t\}$, $t\ge 0$. If $M$ is
further a BMO martingale, then the stochastic exponential $\sE(M)$
is itself a uniformly integrable martingale (see \cite[Theorem
2.3]{Ka}).

The theory of BMO was brought into the study of quadratic BSDEs for
the following reason. Consider, for example, the
BSDE$(T,\xi,g)$ (see (\ref{BSDE})) where the generator $g$ has a
quadratic growth. Assume that
there is some $k>0$ (we may assume without loss of generality that
$k \ge \frac{1}{2}$) such that for $dt \times dP$-a.s. $(t,\o) \in
[0,T] \times \O$,
 \bea
 \label{gquad}
  |g(t,\o,y,z)| \leq k(1+|z|^2), \qq (y,z) \in \hR \times \hR^d,
 \eea
%
and denote $(Y,Z)\in \hC^\infty_\bF([0,T]) \times
\cH^2_\bF([0,T];\hR^d)$ be a solution of the BSDE$(T,\xi,g)$.
For any $\t\in\cM_{0,T}$, applying It\^o's formula to $\dis e^{4kY_t}$ from $\t$ to $T$ one has
 \beas
 e^{4kY_{\t}}+8k^2\int_{\t}^T e^{4kY_s}|Z_s|^2ds &=& e^{4kY_T} +
 4k\int_{\t}^T e^{4kY_s}g(s,Y_s,Z_s)ds-4k\int_{\t}^T e^{4kY_s}Z_sdB_s
 \nonumber \\
 &\leq& e^{4kY_T}+4k^2\int_{\t}^Te^{4kY_s}\big(1+|Z_s|^2\big)ds
  -4k\int_{\t}^T e^{4kY_s}Z_sdB_s.
  \eeas
It is then not hard to derive, using
some standard arguments, the following estimate:
 \bea\label{bmo1}
 E\big[\int_{\t}^T |Z_s|^2 ds|\cF_{\t}\big] \leq e^{4k \|Y \|_\infty}
 E\big[e^{4k\xi}-e^{4kY_{\t}}|\cF_{\t}\big]+e^{8k \| Y \|_\infty}(T-{\t}).
 \eea
In other words, we conclude that
%
$Z\in$ BMO, and that
 \bea
 \label{BMO1}
 \|Z \|^2_{BMO} \leq (1+T)e^{8k \| Y \|_\infty}.
 \eea

\section{Quadratic $g$-Evaluations on $L^\infty(\cF_T)$}
\setcounter{equation}{0}

Our study of the $g$-evaluation/expectation benefited greatly from
the techniques used to treat the quadratic BSDEs, initiated by
Kobylanski \cite{Ko}.
We first list some results regarding the existence, uniqueness, and
comparison theorems for the quadratic BSDEs.
Throughout the rest of the paper we assume that the generator $g$ in
BSDE($T,\xi,g$) (\ref{BSDE}) takes the form:
 \beas
 g(t,\o,y,z)=g_1(t,\o,y,z)y+g_2(t,\o,y,z), \q~~ \fa (t,\o,y,z) \in
 [0,T] \times \O  \times \hR \times \hR^d,
 \eeas
and satisfies the following {\it Standing Assumptions}:
 \bi
 \item[{\bf (H1)}] Both $g_1$ and $g_2$ are $\sP \otimes \sB(\hR) \otimes \sB(\hR^d)$-measurable
 and both $g_1(t,\o,\cdot,\cdot)$ and $g_2(t,\o,\cdot,\cdot)$ are continuous for any $(t,\o) \in
 [0,T] \times \O$;
 \item[{\bf (H2)}] There exist a constant $k>0$ and an increasing function $\ell:
 \hR^+ \mapsto \hR^+$, such that for $dt \times dP$-a.s. $(t,\o) \in
 [0,T] \times \O$, \beas |g_1(t,\o,y,z)| \leq k \q\mbox{and} \q
 |g_2(t,\o,y,z)| \leq k+\ell(|y|)|z|^2,\q~~ (y,z) \in \hR \times
 \hR^d; \eeas
 \item[{\bf (H3)}] With the same increasing function $\ell$, for $dt \times dP$-a.s.
 $(t,\o) \in [0,T] \times \O$,
 \beas
 \Big|\frac{\pa g}{\pa z}(t,\o,y,z)\Big| \leq
\ell(|y|)(1+|z|),\qq (y,z) \in \hR \times \hR^d;
 \eeas
\item[{\bf (H4)}] For any $ \e > 0$, there exists a positive function $h_\e(t) \in L^1[0,T]$
such that for $dt \times dP$-a.s. $(t,\o) \in [0,T] \times \O$,
\beas \frac{\pa g}{\pa y}(t,\o,y,z) \leq h_\e(t)+\e|z|^2,\qq  (y,z)
\in \hR \times \hR^d. \eeas
 \ei

Under the assumptions (H1)-(H4), it is known (cf. \cite[Theorem 2.3
and 2.6]{Ko}) that for any $\xi \in L^\infty(\cF_T)$, the BSDE
(\ref{BSDE}) admits a unique solution
 $(Y,Z) \in \hC^\infty_\bF([0,T]) \times \cH^2_\bF([0,T];\hR^d)$.
In fact, this result can be extended to the following more general
form, which will be useful in our future discussion.

\begin{prop}\label{BSDEV}
Assume that $g$ satisfies (H1)-(H4). For any $\xi \in
L^\infty(\cF_T)$ and any $V \in \hD^\infty_\bF([0,T])$, the BSDE
\bea\label{BSDE1} Y_t=\xi+\int_t^T g(s,Y_s,Z_s)ds +V_T-V_t-\int_t^T
Z_s dB_s, \qq  t \in [0,T], \eea admits a unique solution $(Y,Z) \in
\hD^\infty_\bF([0,T]) \times \cH^2_\bF([0,T];\hR^d)$.
\end{prop}

{\it Proof.} We define a new generator $\tilde{g}$ by
$\tilde{g}(t,\o,y,z) \dfnn g(t,\o,y-V_t(\o),z),~  (t,\o,y,z) \in
[0,T] \times \O \times \hR \times \hR^d$. Then it is easy to see
that for any $(t,\o,y,z) \in [0,T] \times \O \times \hR \times
\hR^d$
 \beas
 \tilde{g}_1(t,\o,y,z)&=& g_1(t,\o,y-V_t(\o),z),\\
 \tilde{g}_2(t,\o,y,z) &=& g_2(t,\o,y-V_t(\o),z)-g_1(t,\o,y-V_t(\o),z)
 V_t(\o).
 \eeas
It can be easily verified that $\tilde{g}$ also satisfies (H1)-(H4).
We can then conclude (see, \cite{Ko}) that the
BSDE$(T,\xi+V_T,\tilde{g})$ admits a unique solution $(\tilde{Y},Z)
\in \hC^\infty_\bF([0,T]) \times \cH^2_\bF([0,T];\hR^d)$. But this
amounts to saying that $(\tilde{Y}-V,Z)$
is the unique solution of (\ref{BSDE1}), proving the corollary.
\qed

\ms

Proposition \ref{BSDEV} indicates that if $g$ satisfies (H1)-(H4),
then we can again define a $g$-evaluation $\cE^g_{s,t}:
L^\infty(\cF_t)\mapsto L^\infty(\cF_s)$ for $0 \le s \le t \le T$,
as in the previous section. We shall name it as the ``{\it quadratic
$g$-evaluation/expectation}" for obvious reasons. More generally,
for any $\si$, $\t \in \cM_{0,T}$ such that $\si \le \t$, $\pas$, we
can define the quadratic $g$-evaluation $\cE^g_{\si,\t}:
L^\infty(\cF_\t) \mapsto L^\infty(\cF_\si)$ by $\cE^g_{\si,\t}[\xi]
\dfnn Y^\xi_\si$, where $\xi \in L^\infty(\cF_\t)$, and $Y^\xi$
satisfies the BSDE:
 \bea\label{BSDEtau}
Y^\xi_t=\xi+\int_t^T \b1_{\{s < \t
\}}g(s,Y^\xi_s,Z^\xi_s)ds-\int_t^T Z^\xi_sdB_s,\qq  t \in [0,T].
 \eea
with $Z^\xi \in \cH^2_\bF([0,T];\hR^d)$,
and
$Y^\xi_t=Y^\xi_{t \land \t}$ and $Z^\xi_t=\b1_{\{t < \t \}}Z^\xi_t$,
$\pas$ In particular, if $\t=T$, we
define the quadratic $g$-expectation of $\xi$ for any $\si \in
\cM_{0,T}$ by $\cE^g[\xi|\cF_\si] \dfnn \cE^g_{\si,T}[\xi]$.

We note that, similar to the deterministic-time case,
$\cE^g_{\si,\t} $ has the following properties: \bi
\item[\(1\)] {\it Time-Consistency:}~~For any $\rho$, $\si$,
$\t \in \cM_{0,T}$ with $\rho \le \si \le \t$, $\pas$, we have
 \beas
\cE^g_{\rho,\si}\big[\cE^g_{\si,\t}[\xi]\big]=\cE^g_{\rho,\t}[\xi],
\q~~ \pas \q \fa \xi \in L^\infty(\cF_\t);
 \eeas

\item[\(2\)] {\it Constant-Preserving:}~~Assume (\ref{g0}),
$\cE^g_{\si,\t}[\xi]=\xi$, $\pas$, $\fa \xi \in
L^\infty(\cF_\si)$;

\item[\(3\)] {\it ``Zero-one Law":}~~For any $\xi \in L^\infty(\cF_\t)$
and $A \in \cF_\si$, we have $\b1_A\cE^g_{\si,\t}[\b1_A \xi]=\b1_A
\cE^g_{\si,\t}[\xi]$, $\pas$; Moreover, if $g(t,0,0)=0$, $\dtp$,
then $\cE^g_{\si,\t}[\b1_A \xi]=\b1_A \cE^g_{\si,\t}[\xi]$, $\pas$;

\item[\(4\)] {\it ``Translation Invariant":}~~If $g$ is independent of $y$, then
 \beas
 \cE^g_{\si,\t}[\xi+\eta]=\cE^g_{\si,\t}[\xi]+\eta,\q \pas \q
\fa \eta \in L^\infty(\cF_\si),\q \xi \in L^\infty(\cF_\t).
 \eeas

\item[\(5\)] {\it Strict Monotonicity:}~~For any $\xi,\eta \in L^\infty(\cF_\t)$ with $\xi \geq \eta$,
$\pas$, we have $\cE^g_{\si,\t}[\xi] \geq \cE^g_{\si,\t}[\eta]$,
$\pas$; Moreover, if $\cE^g_{\si,\t}[\xi] = \cE^g_{\si,\t}[\eta]$,
$\pas$, then $\xi=\eta$, $\pas$
 \ei

We remark that the last property (5) above is not completely
obvious. In fact this will be a consequence of so-called ``strict
comparison theorem" for quadratic BSDEs, a strengthened version of
the usual comparison theorem (see, for example, \cite[Theorem
2.6]{Ko}). For completeness we shall present such a version, under
the following conditions that are similar to those in \cite{Ko}, but
slightly weaker than (H1)--(H4).
 \bi
 \item[{\bf (A1)}]~$g$ is $\sP \otimes \sB(\hR) \otimes \sB(\hR^d)$-measurable
and $g(t,\o,\cdot,\cdot)$ is continuous for any $(t,\o) \in [0,T]
\times \O$;
\item[{\bf (A2)}]~For any $M>0$, there exist $\ell \in L^1 [0,T], k \in L^2 [0,T]$
and $C>0$ such that for $dt \times dP$-a.s. $(t,\o) \in [0,T] \times
\O$ and any $(y,z) \in [-M,M]\times \hR^d$,
 \beas
 \big|g(t,\o,y,z)\big| \leq \ell(t)+C|z|^2 \q\mbox{and}\q
 \Big|\frac{\pa g}{\pa z}(t,\o,y,z)\Big| \leq k(t)+C|z|;
 \eeas
 \item[{\bf (A3)}]~For any $\e>0$, there exists a positive function $h_\e \in L^1[0,T]$
such that for $dt \times dP$-a.s. $(t,\o) \in [0,T] \times \O$ and
any $(y,z) \in \hR \times \hR^d$,
 \beas \frac{\pa g}{\pa y}(t,\o,y,z)\leq h_\e(t)+\e|z|^2.
 \eeas \ei
\begin{thm}
Assume (A1)-(A3). Let $\xi^1, \xi^2 \in L^\infty(\cF_T)$ and $V^i$,
$i=1,2$ be two adapted, integrable, right-continuous processes null at $0$. Let $\big(Y_t^i,Z_t^i\big) \in
\hD^\infty_\bF([0,T]) \times \cH^2_\bF([0,T];\hR^d)$, $i=1,2$
be solutions to the BSDEs:
 \beas
Y^i_t=\xi^i+\int_t^T g(s,Y_s^i,Z_s^i)ds+\int_t^T d V_s^i-\int_t^T
Z_s^i d B_s, \q~~ t \in [0,T],\q i=1,2,
 \eeas
respectively. If $\xi^1 \geq \xi^2$, $\pas$ and $ V_t^1-V_t^2$ is
increasing, then it holds \pas~that
  \bea\label{scomp1}
   Y_t^1 \geq Y_t^2, \qq t \in [0,T].
 \eea
 Moreover, if $Y_\t^1 =
 Y_\t^2$ for some $\t \in \cM_{0,T}$, then it holds \pas ~ that
 \bea\label{scomp2}
 \xi^1 = \xi^2,  \q~~ \mbox{and}\q~~ V^1_T- V^2_T=  V^1_\t - V^2_\t .
 \eea
\end{thm}

  {\it Proof.}   It is not hard to see that (\ref{scomp1}) is a mere generalization of \cite[Theorem 2.6]{Ko}, thus we only need to
prove (\ref{scomp2}). Let $M\dfnn\|Y^1\|_\infty+\| Y^2 \|_\infty$,
and define $\D\eta=\eta^1-\eta^2$ for $\eta=Y$, $Z$, $V$,
respectively. Then $\D Y$ satisfies:
 \bea
 \label{DeltaY}
  d \D Y_t&=&
 -\big(g(t,Y^1_t,Z^1_t)-g(t,Y^2_t,Z^2_t)\big)dt-d \D V_t+\D Z_t d B_t \nonumber\\
 &=& -\int_0^1 \Big(\frac{\pa g}{\pa y}(\Xi_t^\l)\D Y_t + \frac{\pa
 g}{\pa z}(\Xi_t^\l)\D Z_t\Big)  d\l d t -d \D V_t+\D Z_t d B_t\\
 &=& -a_t\D Y_t d t -d \D V_t+\D Z_t(-b_tdt+ d B_t),\nonumber
 \eea
where
$\Xi_t^\l\dfnn (t,\l \D Y_t + Y_t^2 ,\l\D Z_t + Z_t^2)$, and
 \beas
 a_t\dfnn \int_0^1 \frac{\pa g}{\pa y}(\Xi_t^\l) d \l
 \q~\mbox{and}\q~ b_t \dfnn \int_0^1 \frac{\pa g}{\pa z}(\Xi_t^\l) d
 \l, \qq t \in [0,T].
 \eeas
Note that $|\l \D Y_t + Y_t^2|\leq M,~\fa t \in [0,T]$, $\pas$, by
using some standard arguments with the help of assumptions
(A1)--(A3) as well as
%
the Burkholder-Davis-Gundy inequality we deduce from (\ref{DeltaY})
that
 \bea
 \label{bddab}
 E \Big\{ \underset{t\in[0,T]}{\sup}\int_0^t a_sd s + \underset{t\in[0,T]}{\sup} \Big|
 \int_0^t b_sd B_s \Big| \Big\} <
 \infty.
 \eea
Define $Q_t\dfnn \exp \Big\{ \int_0^t a_sd s-\frac{1}{2}\int_0^t
|b_s|^2 d s +\int_0^t b_sd B_s\Big\}$, $t\ge 0$, and
 \beas
 \t_n \dfnn \inf\big\{ t \in [\t,T]: Q_t>n
\big\}\land T, \q~~  n \in \hN,
 \eeas
we see that $\t_n\ua T$, $\pas$, and (\ref{bddab}) indicates that
there exists a null set $\sN$ such that for each $\o\in\sN^c$,
$T=\t_m (\o)$ for some $m \in \hN$. On the other hand, for any $n
\in \hN$, integrating by parts on $[\t,\t_n]$ yields that
 \beas Q_{\t_n}\D Y_{\t_n} &\tneg=\neg\dneg& Q_\t \D Y_\t -
\int_\t^{\t_n} Q_t \D Y_t a_t d t
- \int_\t^{\t_n}  Q_t \D Z_t b_t dt - \int_\t^{\t_n}   Q_t d \D V_t  \\
& &+ \neg \int_\t^{\t_n} \neg   Q_t \D Z_t d B_t
 + \neg \int_\t^{\t_n}\neg \D Y_t Q_t a_t d t +\neg \int_\t^{\t_n}\neg \D Y_t Q_t b_t d B_t
  +\neg \int_\t^{\t_n}\neg   Q_t \D Z_t b_t d t \\
&\tneg=\tneg& - \int_\t^{\t_n}   Q_t d \D V_t + \int_\t^{\t_n}
  Q_t \D Z_t d B_t + \int_\t^{\t_n} \D Y_t Q_t b_t
d B_t.
 \eeas
Taking expectation on both sides gives:
 $$ E\Big\{Q_{\t_n}\D Y_{\t_n} + \int_\t^{\t_n}  Q_t d \D V_t \Big\}=0,
 $$
which implies that there exists a null set $\sN_n$ such that for any
$\o \in {\sN_n}^c$, it holds that $\D Y_{\t_n(\o)}(\o)=0$ and $
\D V_{\t_n(\o)}(\o) = \D V_{\t(\o)}(\o)  $. Therefore, for
any $\o \in \Big\{ \sN \cup \big( \underset{n \in \hN}{\cup}\sN_n
\big)\Big\}^c$, one has
\beas
   \D Y_T(\o)=0 \q \mbox{and}\q \D V_T(\o) = \D V_{\t(\o)}(\o) .
   \eeas
   This completes the proof. \qed

In most of the discussion below, we assume the generator $g$
satisfies (H1)-(H4) (hence (A1)-(A3)). We first extend a property of
$g$-expectations \cite[Proposition 3.1]{BCHMP} to the case of
quadratic $g$-evaluations.

\begin{prop}
Assume (H1)--(H4). Assume further that the generator $g$ is
deterministic. For any $t\in[0,T]$ and $\xi \in L^\infty(\cF_t)$, if $\xi $ is independent of $\cF_s$ for some $s\in[0, t)$, then the random variable $\cE^g_{s,t}[\xi]$
is deterministic.
\end{prop}

{\it Proof:} Let $0\le s<t\le T$ be such that $\xi\in
L^\infty(\cF_t)$ and that it is independent of $\cF_s$. It suffices
to show that $\cE^g_{s,t}[\xi]=c$, $\pas$ for some constant $c$.
To see this, for any $r \in [0,t-s]$, we define $B'_r=B_{s+r}-B_s$,
$\cF'_r=\si\big(B'_u, u\in [0,r]\big)$, and $\bF'=\{\cF'_r\}_{r \in [0,t-s]}$.
Clearly, $B'$ is an $\bF'$-Brownian motion on $[0,t-s]$. Since $\xi
\in \cF_t$ is independent of $\cF_s$, one can easily deduce that
$\xi \in \cF'_{t-s}$. Now we denote by $\{(Y'_r,Z'_r)\}_{r \in
[0,t-s]}$ the unique solution to the BSDE:
 \beas
  Y'_r=\xi+\int_r^{t-s}g(s+u,Y'_u,Z'_u)du-\int_r^{t-s} Z'_u dB'_u,
  \q~~~  r \in [0,t-s].
 \eeas
The simple change of variables $r=v-s$ and $w=s+u$ yields that
 \beas
 Y'_{v-s}&=&\xi+\int_v^t g(w, Y'_{w-s},Z'_{w-s})d w-\int_v^t
 Z'_{w-s}d B'_{w-s} \\
 &=& \xi+\int_v^t g(w, Y'_{w-s},Z'_{w-s})d w-\int_v^t
 Z'_{w-s}d B_w, \q~~~  v \in [s,t].
 \eeas
In other words, $\{(Y'_{v-s},Z'_{v-s})\}_{v \in [s,t]}$ is a
solution to  BSDE$(t,\xi,g)$ on $[s,t]$. The uniqueness of the
solution to BSDE then leads to that
 $Y'_{v-s}=\cE^g_{v,t}[\xi],$ $ v \in [s,t]$. In particular, one has
$\cE^g_{s,t}[\xi]=Y'_0$, $\pas$, which is a constant by the
definition of $\bF'$ and the
 Blumenthal $0$-1 law, completing the proof.
 \qed

 \ms
As we can see from the discussion so far, so long as the
corresponding quadratic BSDE is well-posed, the resulting
$g$-evaluation/expectation should behave very similarly to those
with linear growth generators, with almost identical proofs using
the properties obtained so far. We therefore conclude this section
by listing some further properties of the $g$-evaluation/expection
in one proposition for ready references, and leave the proofs to the
interested reader.

\begin{prop}
\label{gevcomp} Let $g_i$, $i=1,2$, be two generators both satisfy
(H1)-(H4).
\begin{itemize}
\item[{1)}] Suppose that $g_i(t,0,0)=0$, $i=1,2$, and that
 \bea \label{gicomp}
 \cE^{g_1}_{0,t}[\xi]=\cE^{g_2}_{0,t}[\xi],\q~~\fa t \in [0,T],~~~\fa
 \xi \in L^\infty(\cF_t),
 \eea
then for any $ \xi \in L^\infty(\cF_T)$, it holds $\pas$ that
$\cE^{g_1}_{t,T}[\xi]=\cE^{g_2}_{t,T}[\xi]$, $\fa  t \in [0,T]$.


\item[{2)}] Suppose further that $g_i$, $i=1,2$ are independent of $y$, 
For any $t \in [0,T]$, if
$\cE^{g_1}_{0,t}[\xi]\le \cE^{g_2}_{0,t}[\xi]$, $ \fa  \xi \in L^\infty(\cF_t)$, then for any $\xi \in L^\infty(\cF_t)$, it holds
$\pas$ that $\cE^{g_1}_{s,t}[\xi]\le \cE^{g_2}_{s,t}[\xi]$, $\fa s
\in [0,t]$.
\end{itemize}
\end{prop}


\ms

To end this section, we state  a stability result of quadratic BSDEs
which is a slight generalization of Theorem 2.8 in \cite{Ko}. Since there is no
substantial difference in the proof, we omit it.

\begin{thm}
\label{stable} Let $\{g_n\}$ be a sequence of generators satisfying (H1) and (H2) with the
same constant $k>0$ and increasing function $\ell$. Denote, for each $n \in \hN$, $(Y^n,Z^n) \in
\hC^\infty_\bF([0,T]) \times \cH^2_\bF([0,T];\hR^d)$ to be a
solution of BSDE$(T,\xi_n,g_n)$ with $\xi_n \in L^\infty(\cF_T)$.

Suppose that $\{\xi_n\}$ is a bounded sequence in $L^\infty(\cF_T)$, and converges $P$-a.s. to some
$\xi \in L^\infty(\cF_T)$; and that
for $dt \times dP$-a.s. $(t,\o) \in [0,T] \times \O$, $\{g_n(t,\o,y,z)\}$ converges to $g(t,\o,y,z)$ locally uniformly in
$(y,z) \in \hR\times \hR^d$ with $g$ satisfying (H1)-(H4). Then
BSDE$(T,\xi,g)$ admits a unique solution $(Y,Z) \in
\hC^\infty_\bF([0,T]) \times \cH^2_\bF([0,T];\hR^d)$ such that
$P$-a.s. $Y^n_t$ converges to $Y_t$ uniformly in $t \in [0,T]$ and
that $Z^n$ converges to $Z$ in $\cH^2_\bF([0,T];\hR^d)$.
\end{thm}

\section{Some Fine Properties of Quadratic $g$-Evaluations}

In this section we extend some fine properties of $g$-evaluation to
the quadratic case. These properties have been discovered for
different reasons in the linear growth cases, and they form an
integral part of the theory of nonlinear expectation. In the
quadratic case, however, the proofs need to be adjusted, sometimes
significantly. We collect some of them here for the distinguished
importance.

 \ss We begin by a representation theorem for the generators via
quadratic $g$-expectation.

\begin{thm} \label{repre3}
Assume (H1)--(H4).  Let  $ (  t,y,z)  \in   [0,T) \times \hR \times \hR^d$. If $g$ satisfies

 \ss \no \(g1\)   $ \dis \underset{ (s, \,y' ) \to
  (t^+ \neg ,\,y)}{\lim} \, g(s,y',z) =  g(t,y,z)$,  \pas~ and

  \ms   \no \(g2\) For some $\e_0 \in (0,T-t]$ and some $\d>0$, there exists an
integrable process $\{\wt{h}_s\}_{s \in [t,t+\e_0]}$ such that
for $dt \times dP$-a.s. $(s,\o) \in [t,t+\e_0] \times \O$,
 \beas
 \qq  \q   \frac{\pa g}{\pa
y'}(s, y',z) \ge \wt{h}_s,\q \fa y' \in \hR   \hb{ with } |y'-y| \le \d   ,
 \eeas
 then it holds $\pas$ that
 \beas 
  g(t,y,z)=\underset{\e \sa 0}{\lim}\,\frac{1}{\e}\big(\cE^g_{t,(t+\e)\land \t
}[y+z(B_{(t+\e)\land \t}-B_t)]-y\big) ,
 \eeas
where $\t \dfnn\inf\big\{s>t: |B_s-B_t|>\frac{\d}{1+|z|} \big\}\land T$.
\end{thm}

{\it Proof:} 
 We set $M \dfnn 1+ |y|+ \frac{\d |z| }{1+|z|}  $, and $\wt{M}\dfnn k M  + 2 \ell (4M ) |z|^2 $.
 By reducing $\e_0$, we may assume that $\wt{M}  \e_0 e^{k \e_0} \le \frac{\d}{1+|z|} \land  \frac{1}{4\ell (4M ) } $.

  \ms   Fix $\e \in (0,  \frac{\ln 2}{k } \land \e_0]$. Since $\|z(B_{(t+\e)\land
\t}-B_t)\|_\infty \le \frac{\d |z|}{1+|z|}$, there exists a unique solution $\big\{(Y^{\e}_s,Z^{\e}_s)\big\}_{s \in [t,t+\e]} \neg
\in \neg \hC^\infty_\bF([t,t+\e]) \neg \times \neg
\cH^2_\bF([t,t+\e]; \hR^d)$ to the following BSDE:
 \beas
 ~~~ Y^{\e}_s\neg =y\neg +\neg z(B_{(t+\e)\land \t}\neg -B_t)\neg
 +\neg \int_{s }^{t+\e  } \neg \b1_{\{r < \t \}} g(r, Y^{\e}_r, Z^{\e}_r)dr\neg -\neg \int_{s }^{t+\e  } \neg Z^{\e}_r dB_r, ~~~  s \in [t,t+\e].
 \eeas
 We know from Corollary 2.2 of \cite{Ko} that $\|Y^{\e}\|_\infty \le (|y|+ \frac{\d |z|}{1+|z|}+ k\e) e^{k\e}
\le 2 M   $. Now let
 \beas
 \wt{Y}^{\e}_s \dfnn Y^{\e}_s-y-z(B_{s\land \t}-B_t), \q
 \wt{Z}^{\e}_s \dfnn Z^{\e}_s-\b1_{\{s < \t\}} z  , \q \fa s \in [t,t+\e] .
 \eeas
It is easy to check that $\{(\wt{Y}^{\e}_s,\wt{Z}^{\e}_s)\}_{s \in
[t,t+\e]}$ is a solution of the BSDE:
 \bea
 \wt{Y}^{\e}_s=\int_s^{t+\e} \tilde{g}(r,\wt{Y}^{\e}_r,\wt{Z}^{\e}_r)dr- \int_s^{t+\e} \wt{Z}^{\e}_r dB_r,
 \q~~~  s \in [t,t+\e]  \label{BSDErev}
 \eea
  with $ \tilde{g}(s, \o, y',z') \dfnn \vf(y') \b1_{\{s < \t\}}
  g\big(s, \o, y'+y+z(B_{s\land \t}(\o)-B_t(\o) ),z'+z\big)$, $
(s,\o, y',z')\in [t,t+\e]\times \O \times \hR \times \hR^d$ where  $\vf: \hR
\mapsto [0,1]$ is an arbitrary $C^1(\hR)$ function that
  equals to $1$ inside $ [-3M, 3M]$,  vanishes outside
$ (-3M-1 , 3M+1 )$  and satisfies $ \underset{ 3M <|x|<3M+1 }{\sup} | \vf'(x)  |  \le 2 $.
 For any  $(s,\o, y',z')\in [t,t+\e]\times \O \times \hR \times \hR^d$, we see that
 \beas
     \tilde{g}(s,\o, y',z') = \tilde{g}_1(s,\o, y',z') y' +\tilde{g}_2(s,\o, y',z')
 \eeas
with
 \beas
 \tilde{g}_1(s,\o, y',z')   &=  &       \vf(y') \b1_{\{s < \t\}}g_1\big(s, \o, y'+y+z(B_{s\land \t}(\o)-B_t(\o)),z'+z\big),  \\
     \tilde{g}_2(s,\o, y',z') &= &   \vf(y') \b1_{\{s < \t\}} g_1\big(s, \o, y'+y+z(B_{s\land \t}(\o)-B_t(\o)),z'+z\big)  \\
    &&  \times  \big(y+z(B_{s\land \t}(\o)-B_t(\o))\big)  \\
    &&  +   \vf(y')  \b1_{\{s < \t\}} g_2\big(s, \o, y'+y+z(B_{s\land \t}(\o)-B_t(\o)),z'+z\big).
  \eeas
 One can easily deduce from (H2) and (H3) that for $dt \times dP$-a.s. $(s,\o )\in [t,t+\e]\times \O$,  it holds for any $( y',z')\in \hR \times \hR^d$ that
  \bea
   |\tilde{g}_1(s,\o, y',z') |  &\le &    k  \label{estimate1} \\
     |\tilde{g}_2(s,\o, y',z') | 
        &  \le &    k M  + 2 \ell (4M   )\big(|z|^2+ |z'|^2  \big)= \wt{M}+ 2 \ell (4M  )  |z'|^2 \qq  \q     \\
   \hb{and }  \q       \left|  \frac{\pa \tilde{g} }{\pa z'}(s,\o,y',z') \right| 
    &\le&   \ell (4 M  ) \big(1+  |z'|+|z|\big) .   \label{estimate3}
  \eea
    Corollary 2.2 of \cite{Ko} once again shows that $\|\wt{Y}^{\e}\|_\infty \le  \wt{M} \e e^{k\e} \le \wt{M}  \e_0 e^{k \e_0}
\le \frac{\d}{1+|z|} \land \frac{1}{4  \ell (4M   ) }$.  Applying It\^o's formula to $|\wt{Y}^{\e}_s|^2$ we obtain that
 \bea
|\wt{Y}^{\e}_s|^2\neg +\neg \int_s^{t+\e}\neg |\wt{Z}^{\e}_r|^2 dr\neg =
2\neg \int_s^{t+\e}\dneg \wt{Y}^{\e}_r\, \tilde{g}(r,\wt{Y}^{\e}_r,\wt{Z}^{\e}_r)dr
\neg- \neg 2 \neg \int_s^{t+\e} \dneg \wt{Y}^{\e}_r \wt{Z}^{\e}_r dB_r, ~~
 s\in [t,t+\e]. \label{ito1}
 \eea
Using \eqref{estimate1}-\eqref{estimate3} and some standard manipulations one derives easily
that
 \beas
 \qq  && \hspace{-2.5 cm} 2\neg \int_s^{t+\e}\dneg \wt{Y}^{\e}_r\,
\tilde{g}(r,\wt{Y}^{\e}_r,\wt{Z}^{\e}_r)dr
  = \neg 2\neg\int_s^{t+\e}\neg
\wt{Y}^{\e}_r\,\tilde{g}(r,\wt{Y}^{\e}_r,0)dr\neg+ \neg 2\neg \int_s^{t+\e}\neg
\wt{Y}^{\e}_r \wt{Z}^{\e}_r \big(\dneg \int_0^1 \neg \frac{\pa
\tilde{g}}{\pa z'}(r,\wt{Y}^{\e}_r,\l \wt{Z}^{\e}_r)d\l \big) dr\\
&\tneg \le &\tneg 2 \neg \int_s^{t+\e}\neg
|\wt{Y}^{\e}_r|\big( k  |\wt{Y}^{\e}_r|+ \wt{M}    \big)dr\neg  +\neg  2\ell(4M)\neg
\int_s^{t+\e}\neg |\wt{Y}^{\e}_r|\,| \wt{Z}^{\e}_r| \big(1+|z|+ \frac12 |\wt{Z}^{\e}_r|\big)dr\\
&\tneg \le &\tneg  \neg \int_s^{t+\e}\neg
|\wt{Y}^{\e}_r|\big( 2k  |\wt{Y}^{\e}_r|+ 2\wt{M} +  \ell(4M) (1+|z|)^2  \big)dr\neg  +\neg 2  \ell(4M)\neg
\int_s^{t+\e}\neg |\wt{Y}^{\e}_r|\,    |\wt{Z}^{\e}_r|^2 dr\\
&\tneg \le &\tneg C \e^2 + \frac12  \int_s^{t+\e}\neg |\wt{Z}^{\e}_r|^2 dr, \q~~~ s\in
[t,t+\e] ,
 \eeas
where $C$ is a generic constant depending on $ |y|, |z|, \e_0, \d, k$ and $\ell(4M)$,
which may vary from line to line. Taking the conditional
expectation $E[~|\cF_s]$ on both sides of (\ref{ito1}) we have
 \bea
 E\Big\{ \int_s^{t+\e}|\wt{Z}^{\e}_r|^2
dr\big|\cF_s\Big\} \le C\e^2 , \qq s \in [t,t+\e]. \label{zterm}
  \eea

\no Now, taking the conditional expectation in the BSDE (\ref{BSDErev})
we have
 \beas
 \tneg \frac{1}{\e} \wt{Y}^{\e}_t- \tilde{g}(t,0,0)&=& \frac{1}{\e} E\Big\{\neg
\int_t^{t+\e}\dneg
\big(\tilde{g}(r,\wt{Y}^{\e}_r,\wt{Z}^{\e}_r)-\tilde{g}(t,0,0)\big)dr\Big|\cF_t\Big\} \\
& \tneg=&  \frac{1}{\e} E \Big\{\int_t^{t+\e} \Big[\wt{Z}^{\e}_r
\int_0^1 \frac{\pa \tilde{g}}{\pa z'}(
r,\wt{Y}^{\e}_r,\l \wt{Z}^{\e}_r )d\l \\
&& +  \wt{Y}^{\e}_r  \int_0^1  \frac{\pa \tilde{g}}{\pa y'}(r,\l
\wt{Y}^{\e}_r,0)d\l + \tilde{g}(r,0,0)-\tilde{g}(t,0,0)\Big]dr\Big|\cF_t\Big\}.
 \eeas

\no  We know from (g2) and  (H4)  that for $dt \times dP$-a.s. $(s,\o) \in [t,t+\e] \times \O$,
    \beas
   \qq   \wt{h}_s  \le   \frac{\pa g }{\pa y'}   (s, \o, y'+y+z(B_{s\land \t}(\o)-B_t(\o)), z) \le h_1(s)+|z|^2
   \eeas
 holds for any  $ y' \in \hR \hb{ with }  |y'| \le \frac{\d}{1+|z|}$.  It follows that for $dt \times dP$-a.s. $(s,\o) \in [t,t+\e] \times \O$,
   \beas
         \left| \frac{\pa \tilde{g} }{\pa y'}(s,\o,y',0)\right| &\tneg=& \tneg  \left| \vf'(y') \b1_{\{s < \t\}}
  g(s, \o, y'+y+z(B_{s\land \t}(\o)-B_t(\o)), z) \right|  \\
  &\tneg &  \tneg + \left| \vf(y')  \b1_{\{s < \t\}}   \frac{\pa g }{\pa y'}   (s, \o, y'+y+z(B_{s\land \t}(\o)-B_t(\o)), z) \right|  \\
 & &  \hspace{-2.5 cm}   \le  2k(1+4M)+\big(1+2\ell(4M) \big)|z|^2 +|\wt{h}_s| + h_1(s)  \dfnn h_s
     \eeas
     holds for any  $ y' \in \hR \hb{ with }  |y'| \le \frac{\d}{1+|z|}$.
Clearly, $\{h_s\}_{ s \in [t,t+\e_0]}$ is an integrable process. Then  applying \eqref{estimate3}, (\ref{zterm}) and the H\"older Inequality we have
 \bea
  \big|\frac{1}{\e} \wt{Y}^{\e}_t- \tilde{g}(t,0,0)\big|\neg
&\le & \frac{1}{\e} E \Big\{
\int_t^{t+\e}\Big[\ell(4M)\big((1    +|z|) |\wt{Z}^{\e}_r|+\frac12  |\wt{Z}^{\e}_r|^2\big)+
|\wt{Y}^{\e}_r| h_r \Big]dr\Big|\cF_t\Big\} \nonumber \\
&& + E \Big\{\frac{1}{\e}\int_t^{t+\e} \big|\tilde{g}(r,0,0)- \tilde{g}(t,0,0)\big|
dr\Big|\cF_t\Big\}  \nonumber  \\
&\le &  C\big( \e+ \sqrt{ \e}\big)+ \wt{M}   e^{k \e} E \Big[
\int_t^{t+\e}h_r dr\Big|\cF_t\Big]  \nonumber \\
&& + E \Big\{\frac{1}{\e}
\int_t^{t+\e}  \big|\tilde{g}(r,0,0)- \tilde{g}(t,0,0)\big| dr\Big|\cF_t\Big\}.  \label{tilde-g}
 \eea
  As  $\underset{  s \to  t^+ }{\lim}  \b1_{\{s < \t\}}=1$ and  $\underset{  s \to  t^+ }{\lim}  (B_{s\land \t} -B_t ) =0 $, \pas,
  one can deduce from (g1) that
  \beas
    \underset{  s \to  t^+ }{\lim} \, \tilde{g}(s,  0, 0)  =
    \underset{  s \to  t^+ }{\lim} \,    g(s,   y+z(B_{s\land \t} -B_t ), z)
   =   g(t,   y ,  z) =\tilde{g}(t,0,0)   , \q \pas,
  \eeas
  which implies that
   \beas
     \underset{\e \sa 0}{\lim}\, \frac{1}{\e}
\int_t^{t+\e}  \big|\tilde{g}(r,0,0)- \tilde{g}(t,0,0)\big| dr =0, \q \pas
   \eeas
Since $   | \tilde{g}(s,\o,0,0)| \le \wt{M}$  for $dt \times dP$-a.s. $(s,\o) \in [t,t+\e] \times \O$,
   Lebesgue Convergence
Theorem implies that the right hand side of \eqref{tilde-g}  converges
$P$-a.s. to 0 as $\e \to 0^+$. Therefore,
 \beas
 g(t,y,z) &=& \tilde{g}(t,0,0)=\underset{\e \to
0^+}{\lim}\,\frac{1}{\e}\wt{Y}^{\e}_t=\underset{\e \to
0^+}{\lim}\,\frac{1}{\e}(Y^{\e}_t-y) \\
 &=& \underset{\e \sa  0}{\lim}\,\frac{1}{\e}\big(\cE^g_{t,(t+\e)\land \t
}[y+z(B_{(t+\e)\land \t}-B_t)]-y\big),  \q \pas,
 \eeas
 where we used \eqref{BSDEtau} in the last equality. The proof is now complete.    \qed

\if{0}
\begin{cor} Assume that $g$ is independent of $y$.
For any $t \in [0,T)$, if the process $\{g(s,0)\}_{s \in [0,T]}$ is
$\pas$ right-continuous at $t$, then for any $z \in \hR^d$, it holds
$\pas$ that $ \dis g(t,z)=\underset{\e \to
0^+}{\lim}\,\frac{1}{\e}\big(\cE^g[zB_{t+\e}|\cF_t]-zB_t\big)$.
\end{cor}

{\it Proof:} It is easy to check that $Y_s\dfnn \cE^g[zB_{t+\e}
|\cF_s]$, $s \in [t,t+\e]$ solves the BSDE (\ref{BSDErev}) with $
g^z(s,z') \dfnn g(s,z'+z)$, $ \fa (s,z')\in[t,t+\e]\times \hR^d$.
Then a simple application of Lemma \ref{repre} yields the
conclusion. \qed
 \fi


A simple application of the Theorem above gives rise to a
reverse to the Comparison Theorem of quadratic BSDE:

\begin{thm}
Assume that $g_i$, $i=1,2$ satisfy (H1)-(H4) and (\ref{g0}). Let $\, t \in
[0,T)$.  If $\cE^{g_1}[\xi|\cF_t] \le \cE^{g_2}[\xi|\cF_t]$, $\pas$ for any $\xi
\in L^\infty(\cF_T)$, and if both $g_i$ satisfy (g1) and (g2) for any $(y,z) \in \hR \times \hR^d$,  then it holds $\pas$ that
 \beas
 g_1(t,y,z) \le  g_2(t,y,z), \q \fa  (y,z) \in \hR \times \hR^d.
 \eeas
\end{thm}

We also have the following corollary of Theorem \ref{repre3}.
\begin{prop}
Assume that $g$ satisfies  (H1)-(H4) and (\ref{g0}). We also assume that $\pas$, $g(\cdot, y,z)$
is continuous for any $(y,z) \in \hR \times \hR^d$. If
 $g$ satisfies (g1) and (g2) for any $(t, y,z) \in  [0,T) \times \hR \times \hR^d$, then $g$ is independent of $y$ if and only if
 \beas
 \cE^g[\xi+c]=\cE^g[\xi]+c, \q~~ \fa \xi \in
 L ^\infty(\cF_T),~~~ \fa  c \in \hR.
 \eeas
\end{prop}

{\it Proof:} ``$\Ra$": A simply application of Translation
Invariance of quadratic $g$-expectations.

\no ``$\La$": For any $c \in \hR$, we define a new generator
$g^c(t,\o,y,z)\dfnn g(t,\o,y-c,z)$, $\fa (t,\o,y,z) \in [0,T]\times
\O \times \hR \times \hR^d$. It is easy to check that $g^c$
satisfies (H1)-(H4) as well as the other assumptions on $g$ in this
proposition. For any $\xi \in L^\infty(\cF_T)$, let $(Y,Z)$ denote
the unique solution to BSDE$(T,\xi,g)$. Setting $\tilde{Y}_t=Y_t+c$,
$t\in [0,T]$ one obtains that
 \beas
 \tilde{Y}_t=\xi+c+\int_t^T g^c(s,\tilde{Y}_s, Z_s)ds-\int_t^T Z_s
 dB_s, \q~~ \fa t \in [0,T].
 \eeas
Thus, it holds $\pas$ that
 \beas
 \cE^{g^c}[\xi+c|\cF_t]=\tilde{Y}_t =Y_t+c =\cE^g[\xi|\cF_t]+c, \q~~\fa t \in [0,T].
 \eeas
In particular, taking $t=0$ gives that $\cE^{g^c}[\xi]=\cE^g[\xi]$
for any $\xi \in L^\infty(\cF_T)$. Since $g$ satisfies (\ref{g0}),
it easy to see that the condition (\ref{gicomp}) is satisfied for
$g^1\dfnn g$ and $g^2 \dfnn g^c$. Hence, Proposition \ref{gevcomp}
implies that for any $ \xi \in L^\infty(\cF_T)$, it holds $\pas$
that $\cE^g[\xi|\cF_t]=\cE^{g^c}[\xi|\cF_t]$, $\fa t \in [0,T]$.
Applying Theorem \ref{repre3} we see that for any $(t,z) \in
[0,T)\times \hR^d$, it holds $\pas$ that
$g(t,c,z)=g^c(t,c,z)=g(t,0,z)$. Then (H1) implies that for any $t\in
[0,T)$, it holds $\pas$ that $g(t,y,z)=g(t,0,z)$, $\fa (y,z) \in \hR
\times \hR^d$. Eventually, by our assumption, it holds $\pas$ that
$g(t,y,z)=g(t,0,z)$, $\fa (t,y,z) \in [0,T)\times\hR \times \hR^d$.
This proves the proposition. \qed

\medskip
To end this section we extend another important feature of the
$g$-expectation to the quadratic case: The Jensen's Inequality. We
begin by recalling some basic facts for convex functions, and we
refer to Rockafellar \cite{Rock} for all the notions to appear
below.

Recall that
if $F: \hR^n \mapsto \hR$ is a convex function, then  by considering
the convex real function $f(\l) \dfnn F(\l x)-\big(\l
F(x)+(1-\l)F(0)\big)$, $\l \in \hR$, with  $f(0)=f(1)=0$, it is easy
to check that for any $x \in \hR^n$, it holds that
 \bea
 \label{Convex}
 \left\{\ba{lll}
 F(\l x) \le \l F(x)+(1-\l)F(0), \q~~\mbox{if }\, \l \in [0,1],\\
 F(\l x) \ge \l F(x)+(1-\l)F(0), \q~~\mbox{if }\, \l \in (0,1)^c.
 \ea\right.
 \eea

Next, if  $F:\hR \mapsto \hR$ is a convex (real) function, then we
denote by $\pa F$ the {\it subdifferential} of $F$ (see
\cite{Rock}). In particular, for any $x\in
\hR$, $\pa F(x)$ is simply an interval $[F'_{-}(x),F'_{+}(x)]$,
where $F'_-$ and $F'_+$ are left-, and right-derivatives of $F$,
respectively.
The following result is an extension of the linear growth case (cf.
\cite[Proposition 5.2]{BCHMP}).
\begin{thm} Assume that $g$  is independent of $y$ and
 satisfies (H1)-(H4) and (\ref{g0}).  Let $\, t \in [0,T)$.  If
 $g(s, \o, z)$ is convex in $z$ for  $dt \times dP$-a.s. $(s,\o) \in [t,T] \times \O$,   then
  \beas
  F\big(\cE^g[\xi|\cF_t]\big)\le
\cE^g[F(\xi)|\cF_t], \q  \pas
 \eeas
for any $\xi \in L^\infty(\cF_T)$
with $\pa F\big(\cE^g[\xi|\cF_t]\big)\cap (0,1)^c
\ne \es$, $\pas$
\end{thm}

{\it Proof:}   Since both $F'_{-}(x) $ and $F'_{+}(x) $ are non-decreasing functions, we can define another non-decreasing function:
 \beas
    \beta(x)\dfnn \b1_{\{F'_{-}(x)\le
0\}}F'_{-}(x)+\b1_{\{F'_{-}(x)> 0\}}F'_{+}(x) , \q x \in \hR.
 \eeas
Thus $ \beta_t \dfnn \beta\big(\cE^g[\xi|\cF_t]\big)$ is an $\cF_t$-measurable random variable. Since $\beta(x) \in (0,1)^c$ for any  $x\in \hR$ with $\pa F(x) \cap (0,1)^c \ne \es$, it follows that
 \bea   \label{beta_01}
 \beta_t \in (0,1)^c ,  \q  \pas
  \eea
  One can deduce from the convexity of $F$ that
  \bea \label{jensen}
 \beta_t \big(\xi-\cE^g[\xi|\cF_t]\big) \le
 F(\xi)-F\big(\cE^g[\xi|\cF_t]\big).
 \eea
  Since $\xi
\in L^\infty(\cF_T)$, it is clear that $F(\xi)$, $\cE^g[\xi|\cF_t]$,
$F\big(\cE^g[\xi|\cF_t]\big)$ as well as $\beta_t
\big(\xi-\cE^g[\xi|\cF_t]\big)$ are all of $L^\infty(\cF_T)$. Taking $\cE^g[~|\cF_t]$ on both side of (\ref{jensen}), and using
Translation Invariance of quadratic $g$-expectation we have
 \beas
&&\cE^g[\beta_t \xi |\cF_t]-\beta_t \cE^g[\xi|\cF_t] =\cE^g\big[\beta_t
\big(\xi-\cE^g[\xi|\cF_t]\big)\big|\cF_t\big] \\
&\le& \cE^g\big[F(\xi)-F\big(\cE^g[\xi|\cF_t]\big)\big|\cF_t\big]=
\cE^g[F(\xi)|\cF_t]-F\big(\cE^g[\xi|\cF_t]\big),\q \pas
 \eeas
Hence, it suffices to show that $\beta_t \cE^g[\xi|\cF_t] \le
\cE^g[\beta_t \xi |\cF_t]$, $\pas$ To see this, let $Y_t \dfnn
\cE^g[\xi|\cF_t]$, $t \in [0,T]$. As $\beta_t \in \cF_t$, one has
 \beas
 \beta_t Y_s =\beta_t \xi + \int_s^T \beta_t g(r,Z_r)dr- \int_s^T \beta_t
 Z_r dB_r, \q~~ \fa s \in [t,T].
 \eeas
Since $g$ is convex and satisfies (\ref{g0}), using
(\ref{Convex}) and \eqref{beta_01} we obtain
 \beas
 \beta_t Y_s \le \beta_t \xi + \int_s^T  g(r,\beta_t Z_r)dr- \int_s^T \beta_t
 Z_r dB_r=\cE^g[\beta_t \xi|\cF_s], \q~~ \fa s \in [t,T].
 \eeas
 In particular, we have $ \beta_t
\cE^g[\xi|\cF_t] \le \cE^g[\beta_t \xi |\cF_t]$, $\pas$, proving the
Theorem. \qed

\section{Main Results}
\setcounter{equation}{0}

\ms
In this section we prove the main results of this paper regarding
the {\it quadratic $g$-martingales}. To begin with, we give the
following definition. Recall that $\cE^g_{s,t}[\cd]$, $0\le s\le
t\le T$ denotes the $g$-evaluation.
\begin{defn}
An $X \in L^\infty_\bF([0,T])$ is called a ``$g$-submartingale"
\(resp. $g$-supermartingale\) if for any $0 \le s \le t \le T$, it
holds that
 \beas
 \cE^g_{s,t}[X_t] \ge (\mbox{resp. }\le) X_s,  \qq \pas
 \eeas
$X$ is called a $g$-martingale if it is both a $g$-submartingale and
a $g$-supermartingale.
\end{defn}
We should note here
that, in the above the martingale is defined in terms of
quadratic $g$-evaluation, instead of quadratic $g$-expectation as we
have usually seen. This slight relaxation is merely for convenience
in applications. It is clear, however, that if $g$ satisfies
(\ref{g0}), then the quadratic $g$-martingale defined above should
be the same as the one defined via quadratic $g$-expectations,
thanks to (\ref{coin}).

We shall extend three main results for $g$-expectation to the
quadratic case: the Doob-Meyer decomposition, the optional sampling
theorem, and the upcrossing theorem. Although the results look
similar to the existing one in the $g$-expectation literature, the
proofs are more involved, due to the special nature of the quadratic
BSDEs. We shall present these results separately.

\ms


We begin by proving a Doob-Meyer type decomposition theorem for
$g$-martingales.
\begin{thm} \label{gmd}
{\rm (Doob-Meyer Decomposition Theorem)} Assume (H1)--(H4). Let $Y$
be any $g$-submartingale (resp. $g$-supermartingale) that has
right-continuous paths. Then there exist a c\`adl\`ag increasing
(decreasing) process $A$ null at $0$ and a process $Z \in
\cH^2_\bF([0,T];\hR^d)$ such that \beas Y_t=Y_T+\int_t^T
g(s,Y_s,Z_s)ds-A_T+A_t-\int_t^T Z_s dB_s, \qq t\in[0,T]. \eeas
\end{thm}

{\it Proof.} We first assume that $Y$ is a $g$-submartingale. Set $M
\dfnn (  \|Y\|_\infty + kT) e^{kT}$ and $K \dfnn \ell(M+1)$, we let $\f:\hR
\mapsto [0,1]$ be any $C^2(\hR)$ function that equals to $1$ inside
$\big[e^{-2KM},e^{2KM}\big]$ and vanishes outside
$\big(e^{-2K(M+1)},e^{2K(M+1)}\big)$. Let us construct a new
generator: For any $(t,\o,y,z) \in [0,T] \times \O \times \hR \times \hR^d$,
  \beas
 \tilde{g}(t,\o,y,z) \dfnn  \f(y)\Big[2Ky \,g\Big(t,\o,\frac{\ln{(y)}}{2K}, \frac{z}{2Ky}\Big)
 -\frac{|z|^2}{2y}\Big].
  \eeas
One can deduce from (H2) that for $dt \times dP$-a.s. $(t,\o) \in
[0,T] \times \O$,
 \beas \tilde{g}(t,y,z) \leq 2(M+2)kK\f(y)y, \qq
 (y,z) \in \hR \times \hR^d.
 \eeas
\if{0}
$$\f(y) = \left\{
\begin{array}{ll}
1,\q \mbox{if }y \in \big[e^{-2KM},e^{2KM}\big]\\
0,\q \mbox{if }y \notin \big(e^{-2K(M+1)},e^{2K(M+1)}\big)
\end{array}
\right.$$
 \fi
Since $2(M+2)kK\f(y)y$ is Lipschitz continuous in $y$, we can
construct (cf. \cite{Ko}) a decreasing sequence $g_n(t,y,z)$ of
generators uniformly Lipsichitz in $(y,z)$ such that \pas
 \beas g_n(t,y,z) \searrow \tilde{g}(t,y,z), \qq  \fa (t,y,z) \in
 [0,T] \times \hR \times \hR^d .
 \eeas
Now fix $t \in [0,T]$, for any $\xi \in L^\infty(\cF_t)$ with
$\|\xi\|_\infty \leq \|Y\|_\infty$, we define $y_s \dfnn \cE^g_{s,t}[\xi]$, $s
\in [0,t]$. It follows from \cite[Corollary 2.2]{ Ko} that $\|y\|_\infty \le (  \|Y\|_\infty + kT) e^{kT}=M$.  Applying It\^o's formula we see that $\tilde{y}_s \dfnn
e^{2Ky_s}$, $s \in [0,t]$ together with a process $\tilde{z} \in
\cH^2_\bF([0,t];\hR^d)$ is a solution of the following BSDE:
 \beas \tilde{y}_s=e^{2K\xi}+\int_s^t
 \tilde{g}(r,\tilde{y}_r,\tilde{z}_r)dr-\int_s^t \tilde{z}_r dB_r,
 \qq \fa s \in [0,t].
 \eeas
Since $g_n$ is Lipschitz, a standard comparison theorem implies that
 \beas e^{2K \cE^g_{s,t}[\xi]}=\tilde{y}_s \leq
 \cE^{g_n}_{s,t}[e^{2K\xi}], \q~~~ s \in [0,t], \q~ \pas
 \eeas
In particular, taking $\xi =Y_t$ shows that
 \beas
 e^{2K Y_s} \le e^{2K \cE^g_{s,t}[Y_t]} \le \cE^{g_n}_{s,t}[e^{2K Y_t}],\q s
 \in [0,t], \q~ \pas
 \eeas
Namely, $\tilde{Y} = e^{2KY}$ is a right-continuous
$g_n$-submartingale in the sense of $g^n$-evaluation for any $n \in
\hN$. Applying the known $g$-submartingale decomposition theorem for
the Lipschitz case (see \cite[Theorem 3.9]{Pln}), we can find a
c\`adl\`ag increasing process $A^n$ null at $0$ and a process $Z^n
\in \cH^2_\bF([0,T];\hR^d)$ such that
 \bea
 \label{gnrep}
 \tilde{Y}_t=\tilde{Y}_T+\int_t^T g_n(s,\tilde{Y}_s,Z^n_s)ds
 -A^n_T+A^n_t-\int_t^T Z^n_s dB_s, \qq t\in[0,T],
 \eea
from which we see that $\tilde{Y}$, whence $Y$ is c\`adl\`ag. Note
that, in the representation (\ref{gnrep}), the martingale parts must
coincide for any $m$ and $n$. In other words, one must have
$Z^m=Z^n$ as the elements in $\cH^2_\bF([0,T];\hR^d)$. Thus, for any
$n \in \hN$, (\ref{gnrep}) can be rewritten as
 \beas
 \tilde{Y}_t=\tilde{Y}_T+\int_t^T g_n(s,\tilde{Y}_s,\tilde{Z}_s)ds
 -A^n_T+A^n_t-\int_t^T \tilde{Z}_s dB_s, \qq
 t\in[0,T].
 \eeas
Since $g_n\searrow \tilde{g}$, the Lebesgue Convergence Theorem
implies that
 \beas
 \int_0^T \big[g_n(s,\tilde{Y}_s,\tilde{Z}_s)-\tilde{g}(s,\tilde{Y}_s,
 \tilde{Z}_s)\big]ds \to 0, \q~~ \pas
 \eeas
Consequently, it holds $\pas$ that
 \beas
 A^n_t \to \tilde{A}_t \dfnn \tilde{Y}_t-\tilde{Y}_0 +\int_0^t
 \tilde{g}(s,\tilde{Y}_s,\tilde{Z}_s)ds-\int_0^t
 \tilde{Z}_s dB_s, \q~~ \fa t \in [0,T].
 \eeas
It is easy to check that  $\tilde{A}$ is also a c\`adl\`ag increasing
process null at $0$.
Now let us define a new $C^2(\hR)$
function $\psi$ by   $ \dis  \psi(y) \dfnn   \frac{\f(y)\ln{(y)}}{2K}$, $y \in \hR$.
Applying It\^o's formula to $\psi(\tilde Y_t)$ from $t$ to $T$ one
has
 \beas
 Y_t&  =  & Y_T+  \int_{t+}^T\frac{1}{2K\tilde{Y}_{s-}}\big[\tilde{g}(s,\tilde{Y}_s,\tilde{Z}_s)
 ds-d\tilde{A}_s-\neg \tilde{Z}_sdB_s\big] +\frac{1}{2}\int_{t+}^T
 \frac{|\tilde{Z}_s|^2}{2K\tilde{Y}^2_{s-}}ds\\
 && - \underset{s\in (t,T]}{\sum} \big\{\Delta Y_s - \frac{\Delta
 \tilde{Y}_s}{2K\tilde{Y}_{s-}}\big\} \\
 &=  &Y_T+\neg \int_t^T
 \frac{1}{2K\tilde{Y}_s}\big[\tilde{g}(s,\tilde{Y}_s,\tilde{Z}_s)ds-d\tilde{A}^c_s-\tilde{Z}_sdB_s\big]
 +\frac{1}{2}\int_t^T \frac{|\tilde{Z}_s|^2}{2K\tilde{Y}^2_s}ds-\dneg\underset{s \in (t,T]}{\sum}\Delta Y_s \\
 & = &Y_T+ \neg \int_t^T g(s, Y_s,
 \frac{\tilde{Z}_s}{2K\tilde{Y}_s})ds-\int_t^T\frac{1}{2K\tilde{Y}_s}d\tilde{A}^c_s
 -\int_t^T\frac{\tilde{Z}_s}{2K\tilde{Y}_s}dB_s- \underset{s \in (t,T]}{\sum}\Delta Y_s,
 \eeas
where the second equality is due to the fact that $\Delta
\tilde{Y}_s= \Delta \tilde{A}_s>0$ and $\tilde{A}^c$ denotes the
continuous part of $\tilde{A}$. Clearly, $ A_t \dfnn
\int_0^t\frac{1}{2K\tilde{Y}_s}d\tilde{A}^c_s + \tneg \underset{s
\in (0,t]}{\sum}\Delta Y_s $ is a c\`adl\`ag increasing process null
at $0$, finally we get \beas Y_t=Y_T+\int_t^T g(s, Y_s,
Z_s)ds-A_T+A_t-\int_t^T Z_s dB_s, \qq t\in[0,T]. \eeas

\ms

On the other hand, if $Y$ is a $g$-supermartingale, then one can
easily check that $-Y$ is correspondingly a $g^{-}$-submartingale
with \bea\label{gneg} g^{-}(t,\o,y,z)\dfnn -g(t,\o,-y,-z),\qq \fa
(t,\o,y,z) \in [0,T] \times \O \times \hR \times \hR^d. \eea
Clearly, $g^{-}$ also satisfies (H1)-(H4), thus there exist a
c\`adl\`ag increasing process $A$ null at $0$ and a process $Z \in
\cH^2_\bF([0,T];\hR^d)$ such that \beas
  -Y_t=-Y_T+\int_t^T g^{-}(s,-Y_s,Z_s)ds-A_T+A_t-\int_t^T Z_s dB_s, \qq  t\in[0,T].
\eeas We can rewrite this BSDE as:  \beas Y_t=Y_T+\int_t^T
g(s,Y_s,-Z_s)ds-(-A_T)+(-A_t)-\int_t^T (-Z_s) dB_s, \qq  t\in[0,T].
\eeas The proof is now complete. \qed

\if{0}
\ms The following corollary is significant.
\begin{cor}
Assume that $g$ satisfies (\ref{g0}) and satisfies
(H2) with $\ell$ being constant in (H2). For any right-continuous $Y
\in L^\infty_\bF([0,T])$ and any $\z \in L^\infty_\bF([0,T];\hR^d)$,
if the process $Y_t+\int_0^t \z_s dB_s$, $t \in [0,T]$ is a
$g$-submartingale \(resp. $g$-supermartingale\) in the sense of
$g$-expectation, then there exist a c\`adl\`ag increasing \(resp.
decreasing\) process $A$ null at $0$ and a process $Z \in
\cH^2_\bF([0,T];\hR^d)$ such that \beas Y_t=Y_T+\int_t^T
g(s,Y_s,Z_s+\z_s)ds-A_T+A_t-\int_t^T Z_s dB_s, \qq t\in[0,T]. \eeas
\end{cor}

{\it Proof.} Since $g$ satisfies (\ref{g0}), $Y_\cdot+\int_0^\cdot
\z_s dB_s$, is also a $g$-submartingale \(resp.
$g$-supermartingale\) in the sense of quadratic $g$-evaluation. Then
a BSDE transformation shows that $Y$ is correspondingly a
$\tilde{g}$-submartingale \(resp. $\tilde{g}$-supermartingale\) in
the sense of $\tilde{g}$-evaluation with \beas
\tilde{g}(t,\o,y,z)\dfnn g(t,\o,y,z+\z_t(\o)),\qq \fa (t,\o,y,z) \in
[0,T] \times \O \times \hR \times \hR^d. \eeas It is easy to check
that $\tilde{g}$ satisfies (H1)-(H4). Hence Theorem \ref{gmd}
implies the conclusion. \qed
\fi

We now turn our attention to the {\it Optional Sampling Theorem}. We
begin by presenting a lemma that will play an important role in the
proof of the Optional Sampling Theorem.
\begin{lem}
Let $\t \in \cM_{0,T}$ be
finite valued in a set $0 = t_0 < t_1< \cdots <t_n= T$. If $t_i \le
s< t \le t_{i+1}$ for some $i \in \{0,1, \cdots n-1\}$, then for any $ \xi \in \cF_{t \land \t}$
\bea\label{op0} \cE^g_{s \land \t , t \land \t}[\xi]=\b1_{\{\t \le
t_i\}}\xi+\b1_{\{\t \ge t_{i+1}\}}\cE^g_{s, t}[\xi],\q~~ \pas
 \eea
\end{lem}

{\it Proof.} For any $\xi \in \cF_{t \land \t}$, let $(Y,Z)$ be the unique solution
to the BSDE (\ref{BSDEtau}) with $\t = t \land \t$.
Then we have
 \beas
 \cE^g_{r \land \t, t \land \t}[\xi]=Y_{r \land \t} &=& \xi + \int_{r \land \t}^T \b1_{\{u <  t \land \t \}}
 g(u,Y_u,Z_u)du-\int_{r \land \t}^T \b1_{\{u <  t \land \t \}} Z_udB_u\\
& =& \xi  + \int_r^t \b1_{\{u <  \t \}}g(u,Y_{u \land
 \t},Z_u)du-\int_r^t \b1_{\{u <  \t \}}Z_udB_u, \qq \fa r \in [0,t].
 \eeas
For any $r \in [s,t]$, since $\{ \t \le t_i \}=\{ \t \ge t_{i+1}
\}^c \in \cF_{t_i} \subset \cF_r$,
one can deduce that
 \bea
 \label{op1}
 \b1_{\{ \t \le t_i \}}Y_{r \land \t}\neg&\neg =\neg &\neg
 \b1_{\{ \t \le t_i \}}\xi  + \int_r^t \b1_{\{ \t \le
 t_i \}}\b1_{\{u <  \t \}}g(u,Y_{u \land \t},
 Z_u)du-\int_r^t \b1_{\{ \t \le t_i \}} \b1_{\{u < \t
 \}}Z_udB_u \nonumber\\
 \neg&\neg=\neg&\neg \b1_{\{ \t \le t_i \}}\xi,
 \eea
and that
 \bea
 \label{part1}
 \b1_{\{ \t \ge t_{i+1} \}}Y_{r \land \t}
 &\dneg=\dneg&\b1_{\{ \t \ge t_{i+1} \}}\xi  + \int_r^t \b1_{\{ \t
 \ge t_{i+1} \}} \b1_{\{u <  \t \}}g(u,Y_{u \land \t},Z_u)du-\int_r^t
 \b1_{\{ \t \ge t_{i+1} \}}
 \b1_{\{u <  \t \}}Z_udB_u  \nonumber \\
 &\dneg=\dneg&\b1_{\{ \t \ge t_{i+1} \}}\xi  + \int_r^t \b1_{\{ \t
 \ge t_{i+1} \}}g(u,Y_{u \land \t},Z_u)du-\int_r^t \b1_{\{ \t \ge
 t_{i+1} \}} Z_udB_u.
 \eea
On the other hand, we let $Y'_r=\cE^g_{r,t}[\xi],~r \in [0,t]$. Then
for any $r \in [s,t]$, by the definition of quadratic
$g$-evaluation, one has
 \bea
 \label{part2} \b1_{\{\t \le t_i\}}Y'_r =
 \b1_{\{\t \le t_i\}}\xi + \int_r^t \b1_{\{\t \le
 t_i\}}g(u,Y'_u,Z'_u)du-\int_r^t \b1_{\{\t \le t_i\}}Z'_udB_u.
 \eea
Adding (\ref{part2}) to (\ref{part1}) shows that $\tilde{Y}_r \dfnn
\b1_{\{\t \ge t_{i+1}\}}Y_{r \land \t}+\b1_{\{\t \le t_i\}}Y'_r$ and
$\tilde{Z}_r \dfnn \b1_{\{ \t \ge t_{i+1} \}} Z_r+\b1_{\{\t \le
t_i\}}Z'_r$ solve the following BSDE
 \beas \tilde{Y}_r= \xi+
 \int_r^t g(u,\tilde{Y}_u,\tilde{Z}_u)du-\int_r^t \tilde{Z}_u dB_u,
 \qq \fa r \in [s,t].
 \eeas
Then it is not hard to check that $\hat{Y}_r=\b1_{\{r\ge
s\}}\tilde{Y}_r+\b1_{\{r< s\}}\cE^g_{r,s}[\tilde{Y}_s],~r \in [0,t]$
is the unique solution of BSDE$(t,\xi,g)$. Hence we can rewrite
$\hat{Y}_r=\cE^g_{r,t}[\xi],~ r \in [0,t]$. In particular, it holds
$P$-a.s. that
 \bea
 \label{part3} \b1_{\{\t \ge t_{i+1}\}}Y_{s
 \land \t}=\b1_{\{\t \ge t_{i+1}\}}\tilde{Y}_s=\b1_{\{\t \ge
 t_{i+1}\}}\hat{Y}_s=\b1_{\{\t \ge t_{i+1}\}}\cE^g_{s,t}[\xi].
 \eea
Letting $r=s$ in (\ref{op1}) and then adding it to (\ref{part3}),
the lemma follows. \qed

We are now ready to prove the Optional Sampling Theorem.
\begin{thm}
Assume (H1)-(H4). For any $g$-submartingale $X$ \(resp.
$g$-supermartingale, $g$-martingale\) such that $\underset{\o \in
\O} {\esssup}\underset{t \in [0,T]}{\sup} |X(t,\o)|<\infty$, and for
any $\si$, $\t \in \cM_{0,T}$ with $\si \le \t$, $\pas$ Assume either that
  $\si$ and $\t$ are finitely valued or that $X$ is
right-continuous, then \beas \cE^g_{\si,\t}[X_\t]\geq (\mbox{resp.
}\leq,\,= )\,X_\si, \q~~~ \pas \eeas
\end{thm}

{\it Proof.} We shall consider only the $g$-submartingale case, as
the other cases can be deduced easily by standard argument. To begin
with, we assume that $\t$ takes values in a finite set $0 = t_0 <
t_1< \cdots <t_n = T$. Note that if $t \geq t_n$, then it is clear
that $\cE^g_{t \land \t,\t}[X_\t]=\cE^g_{\t,\t}[X_\t]=X_\t$, $\pas$
We can then argue inductively that for any $t \in [0,T]$,
 \bea
 \label{op2}
 \cE^g_{t \land \t,\t}[X_\t] \ge X_{t \land \t},
 \q~~~ \pas
 \eea
In fact, assume that for some $i \in \{1,\cdot \cdot\cdot n\}$,
(\ref{op2}) holds for any $t \geq t_i$. Then for any $t \in
[t_{i-1},t_i)$, the time-consistence and the monotonicity of
quadratic $g$-evaluations as well as (\ref{op0}) imply that
 \beas
 \cE^g_{t \land \t,\t}[X_\t]&=&
 \cE^g_{t \land \t,t_i \land \t}\big[\cE^g_{t_i \land \t,\t}[X_\t]\big]
 \ge  \cE^g_{t \land \t,t_i \land \t}[X_{t_i \land \t}]\\
 &=& \b1_{\{\t \le t_{i-1}\}}X_{t_i \land \t}
 +\b1_{\{\t \ge t_i\}}\cE^g_{t, t_i}[X_{t_i \land \t}]\\
 &=& \b1_{\{\t \le t_{i-1}\}}X_{t \land \t}+\b1_{\{\t \ge
 t_i\}}\cE^g_{t, t_i}[X_{t_i \land \t}], \qq \pas
 \eeas
Since $\{\t \ge t_i\}=\{\t \le t_{i-1}\}^c \in \cF_t$, the
``zero-one law" of quadratic $g$-evaluations shows that
 \beas
 \b1_{\{\t \ge t_i\}}\cE^g_{t, t_i}[X_{t_i \land \t}] &=&\b1_{\{\t
 \ge t_i\}}\cE^g_{t, t_i}[\b1_{\{\t \ge t_i\}}X_{t_i \land \t}]
 =\b1_{\{\t \ge t_i\}}\cE^g_{t, t_i}[\b1_{\{\t \ge t_i\}}X_{t_i}]\\
 &=&\b1_{\{\t \ge t_i\}}\cE^g_{t, t_i}[X_{t_i}]\ge \b1_{\{\t \ge
 t_i\}}X_t=\b1_{\{\t \ge t_i\}}X_{t \land \t},\qq \pas
 \eeas
Hence, (\ref{op2}) holds for any $t \geq t_{i-1}$, this completes
the inductive step. If $\si$ is also finitely valued, for example in
the set $0 = s_0 < s_1< \cdots <s_m = T$, then it holds $P$-a.s.
 \bea
 \label{op4}
 \cE^g_{\si,\t}[X_\t]=\cE^g_{\si\land\t,\t}[X_\t]=\sum^m_{j=0}\b1_{\{\si
 = s_j\}}\cE^g_{s_j\land \t,\t}[X_\t]\ge \sum^m_{j=0}\b1_{\{\si = s_j\}}X_{s_j\land
 \t}=X_{\si \land \t}=X_\si.
 \eea

For a general $\t \in \cM_{0,T}$, we define two sequences
$\{\si_n\}$ and $\{\t_n\}$ of finite valued stopping times such that
$P$-a.s.
 \beas
 \si_n \searrow \si, \q~~ \t_n \searrow \t, \q ~~ \mbox{and}
 \q~~ \si_n \le \t_n, \q \fa n \in \hN.
 \eeas
Fix $n \in \hN$ and let $(Y^n,Z^n)$ be the unique solution
to the BSDE (\ref{BSDEtau}) with $\xi=X_{\t_n}  $ and $\t =  \t_n$. We know from (\ref{op4}) that $P$-a.s.
 \beas
  Y^n_{\si_m}= \cE^g_{\si_m,\t_n}[X_{\t_n}] \ge X_{\si_m}, \qq \fa m \ge n.
 \eeas
In light of the right-continuity of $X$ and $Y^n$, letting $m \to \infty$ gives that
 \beas
  Y^n_\si    \ge X_\si, \qq \pas
  \eeas
Now let $(Y,Z)$ be the unique solution
to the BSDE (\ref{BSDEtau}) with $\xi=X_{\t}  $. It is easy to see that for $dt \times dP$-a.s. $(t,\o) \in [0,T]
\times \O,~\b1_{\{t \le \t_n\}}g(t,\o,y,z)$ converges to $\b1_{\{t
\le \t\}}g(t,\o,y,z)$ uniformly in $(y,z) \in \hR\times
\hR^d$.  Theorem \ref{stable} then implies that $P$-a.s.
$Y^n_t $ converges to $  Y_t $ uniformly in $t
\in [0,T]$. Thus we have
 \beas
    \cE^g_{\si,\t}[X_\t] = Y_\si   =  \underset{n \to \infty}{\lim} Y^n_\si   \ge  X_\si, \qq \pas,
 \eeas
 proving the theorem.
\qed

Finally, we study the so-called {\it Upcrossing Inequality} for
quadratic $g$-submartingales, which would be essential
for the study of path regularity of $g$-submartingales.

\begin{thm}\label{ucism}
Given a $g$-submartingale $X$,
we set $J \dfnn \big(\|X\|_\infty+kT\big)e^{kT}$ and denote $\wt{X}_t
 = X_t+k(J+1)t,~t \in [0,T]$. As usual, for any finite
set $\cD=\{0 \le t_0<t_1<...<t_n \le T\}$, we let
$U^b_a(\wt{X},\cD)$ denote the number of upcrossings of the interval
$[a,b]$ by $\wt{X}$ over $\cD$. Then there is a BMO
process $\big\{\beta_\cD(t) \big\}_{t \in [0,t_n]} $ such that
 \beas
 E\Big[U^b_a(\wt{X},\cD)\exp{\big(\int_0^{t_n} \beta_\cD(s) dB_s -
 \frac{1}{2} \int_0^{t_n} |\beta_\cD(s)|^2 ds \big)} \Big] \le
 \frac{\|X\|_\infty+k(J+1)T+|a|}{b-a},
 \eeas
and that
 $E\int_0^{t_n} |\beta_\cD(s)|^2 ds \le C$, a constant independent of the choice of $\cD$.
\end{thm}

{\it Proof.} For any $j \in \{1,\cdots n\}$ we consider the
following BSDE:
 \beas
 Y^j_t=X_{t_j}+\int_t^{t_j}
 g(s,Y^j_s,Z^j_s)ds-\int_t^{t_j} Z^j_sdB_s, \qq \fa t \in [t_{j-1},t_j].
 \eeas
Applying Corollary 2.2 of \cite{Ko} one has
 \bea
 \label{usis1}
 \|Y^j\|_\infty \le
 \big(\|X_{t_j}\|_\infty+k(t_j-t_{j-1})\big)e^{k(t_j-t_{j-1})}\le J.
 \eea

Now let us define a $d$-dimensional process
$\beta_\cD(t)=(\beta^1_t,\cdots \beta^d_t)$, $t \in [0,t_n]$ by
 \beas
 \beta^l_t \dfnn \sum^n_{j=1}\b1_{t \in (t_{j-1},t_j]} \int_0^1
 \frac{\pa g}{\pa z_l} \big(t,Y^j_t,(Z^{j,1}_t, \cdots
 \l Z^{j,l}_t ,0,\cdots0) \big) d\l,\qq l \in \{1,\cdots,d\}.
 \eeas
It is easy to see from Mean Value Theorem that for any $t \in
(t_{j-1},t_j]$,
 \bea
 \label{uci2}
 && g(t,Y^j_t,Z^j_t)-g(t,Y^j_t,0)\nonumber\\
 &=&\sum^d_{l=1}\Big\{g\big(t,Y^j_t,(Z^{j,1}_t, \cdots
 Z^{j,l}_t,0,\cdots, 0)\big)\neg
 -\neg g\big(t,Y^j_t,(Z^{j,1}_t, \cdots  Z^{j,l-1}_t,0,\cdots,0)\big) \Big\} \nonumber\\
 &=&\sum^d_{l=1} Z^{j,l}_t \beta^l_t =\lan Z^j_t, \beta_\cD(t)\ran.
 \eea
Moreover,  (H3) implies that
 \bea
 \label{uci4}
 \big|\beta^l_t\big|\le \ell(J)\sum^n_{j=1}\b1_{t \in (t_{j-1},t_j]}(1+ |Z^j_t|),\q~~
  t \in [0,t_n],\q l \in \{1,\cdots,d\}.
 \eea
  We see from  \eqref{BMO1}  that each $Z^j$ is a BMO
process, thus so is $\beta_\cD$. In fact, for any $\t \in
\cM_{0,t_n}$, one can deduce from (\ref{uci4}) that \if{0}
 \bea\label{uci3} E\big[
 \int_\t^{t_n} |\beta_\cD(s)|^2 ds \big| \cF_\t \big]
 &\tneg \le \tneg& 2d \ell(J)^2t_n +2d \ell(J)^2\sum_{j=1}^n
 E\big[\int_{(\t \vee t_{j-1})\land t_j }^{t_j} \nts \nts |Z^j_s|^2 ds  \big|  \cF_\t \big]\\
 &\tneg \leq \tneg& 2d \ell(J)^2T+2d \ell(J)^2\sum_{j=1}^n \Big\{\b1_{\{\t \le t_{j-1}\}}
 E\big[\int_{t_{j-1}}^{t_j} \nts \nts |Z^j_s|^2 ds  \big|  \cF_{\t \land t_{j-1}}\big]\nonumber\\
 & &+\b1_{\{t_{j-1} < \t \leq t_j\}}E\big[\int_{(\t \vee t_{j-1})\land t_j}^{t_j}
 \nts \nts |Z^j_s|^2 ds  \big|  \cF_{(\t \vee t_{j-1})\land t_j}\big]\Big\} \nonumber\\
 &\tneg \leq \tneg& 2d \ell(J)^2 T+2d \ell(J)^2\sum_{j=1}^n \Big\{\b1_{\{\t \le t_{j-1}\}}
 E\big[E[\int_{t_{j-1}}^{t_j} \nts \nts |Z^j_s|^2 ds|\cF_{t_{j-1}}]
 \big|  \cF_{\t \land t_{j-1}}\big]\nonumber\\
 &&+\b1_{\{t_{j-1} < \t \le t_j \}}\norm Z^j_s \norm^2_{BMO_2} \Big\}\nonumber\\
 &\tneg \leq \tneg& 2d \ell(J)^2T +2d \ell(J)^2\sum_{j=1}^n \norm
 Z^j_s \norm^2_{BMO_2} \nonumber
 \eea
 \fi
 \bea
 \label{uci3} && E\big[
 \int_\t^{t_n} |\beta_\cD(s)|^2 ds \big| \cF_\t \big]
 \le Ct_n +C\sum_{j=1}^n E\big[\int_{(\t \vee t_{j-1})\land t_j }^{t_j}
 \nts \nts |Z^j_s|^2 ds  \big|  \cF_\t \big]\\
 &\tneg \leq \tneg& CT\neg+\neg C\sum_{j=1}^n \Big\{\b1_{\{\t \le t_{j-1}\}}
 E\big[\neg \int_{t_{j-1}}^{t_j} \neg \nts \nts |Z^j_s|^2 ds  \big|  \cF_{\t \land t_{j-1}}\big]
 \neg+\neg \b1_{\{t_{j-1} < \t \leq t_j\}}
 E\big[\neg \int_{(\t \vee t_{j-1})\land t_j}^{t_j} \nts \nts |Z^j_s|^2 ds
 \big|  \cF_{(\t \vee t_{j-1})\land t_j}\big]\Big\} \nonumber\\
 &\tneg \leq \tneg& CT\neg+\neg C\sum_{j=1}^n \Big\{\b1_{\{\t \le t_{j-1}\}}
 E\big[E[\int_{t_{j-1}}^{t_j} \nts \nts |Z^j_s|^2 ds|\cF_{t_{j-1}}]  \big|
 \cF_{\t \land t_{j-1}}\big]\neg+\neg \b1_{\{t_{j-1} < \t \le t_j \}}
 \norm Z^j_s \norm^2_{BMO_2} \Big\}\nonumber\\
 &\tneg \leq \tneg& CT \neg+\neg C\sum_{j=1}^n \norm Z^j_s
 \norm^2_{BMO_2}, \nonumber
 \eea
where $C\dfnn 2d \, \ell(J)^2$. Thus $\big\{\sE\big( \beta_\cD \bullet
B \big)_t\big\}_{t \in [0,t_n]}$ is a uniformly integrable martingale.
By Girsanov's theorem we can find an equivalent probability $Q^\cD$
such that
 \if{0}
 \beas \frac{dQ^\cD}{dP}=\exp{\Big\{
 \int_0^{t_n}\beta_\cD(s)dB_s-\frac{1}{2}
 \int_0^{t_n}|\beta_\cD(s)|^2ds \Big\}} \eeas
 \fi
$dQ^\cD/dP=\sE\big(\beta_\cD \bullet B \big)_{t_n}$.
 \if{0}
 Moreover, it is well-known that $\dis B^\cD_t \dfnn B_t-\int_0^t
 \beta_\cD(s) ds~ t \in [0,t_n]$ is the Brownian Motion under
 $Q^\cD$.
 \fi
Then (\ref{uci2}) and (H2) show that for any $j \in\{1,...,n\}$ and
any $t \in [t_{j-1},t_j]$,
 \beas
 Y^j_t &=&X_{t_j} + \int_t^{t_j}\big[g(s,Y^j_s,0)+\lan Z^j_s, \beta_\cD(s)\ran\big]ds-\int_t^{t_j} Z^j_sdB_s\\
 &=& X_{t_j}+ \int_t^{t_j} g(s,Y^j_s,0)ds-\int_t^{t_j} Z^j_s dB^\cD_s \\
 &\le& X_{t_j}+ k(J+1)(t_j-t)-\int_t^{t_j} Z^j_s dB^\cD_s,
 \eeas
where $B^\cD$ denotes the Brownian Motion under $Q^\cD$. Taking the
conditional expectation $E_{Q^\cD}[\cdot|\cF_t]$ on both sides of
the above inequality one can obtain that
 \beas
 \cE^g_{t,t_j}[X_{t_j}]=Y^j_t \le
 E_{Q^\cD}[X_{t_j}|\cF_t]+k(J+1)(t_j-t), \q~~\pas \q \fa t \in
 [t_{j-1},t_j].
 \eeas
In particularly, taking $t=t_{j-1}$ we have
 \beas
 X_{t_{j-1}} \le \cE^g_{t_{j-1},t_j}[X_{t_j}] \le
 E_{Q^\cD}[X_{t_j}|\cF_{t_{j-1}}]+k(J+1)(t_j-t_{j-1}),\q~~ \pas
 \eeas
Hence $\{\wt{X}_{t_j}\}^n_{j=0}$ is a ${Q^\cD}$-submartingale.
Applying the classical upcrossing theorem one has
 \beas
 E_{Q^\cD}\big[U^b_a(\wt{X},\cD)\big] \leq
 \frac{E_{Q^\cD}\big[(\wt{X}_{t_n}-a)^+\big]}{b-a} \le
 \frac{\|X\|_\infty+k(J+1)T+|a|}{b-a}
 \eeas
Furthermore, we denote $C>0$ to be a generic constant depending only
on $d,T, J, k, \|X\|_\infty$, and is allowed to vary from line to line.
Letting $\t=0$ in (\ref{uci3}) one can deduce that
 \beas
  &&   E\int_0^{t_n} |\beta_\cD(s)|^2 ds  \le  C +C\sum_{j=1}^n E\int_{t_{j-1}}^{t_j} |Z^j_s|^2 ds \\
 &\le& C +C\sum_{j=1}^n \Big\{ e^{4\tilde{K}J} E\big[e^{4\tilde{K}Y^j_{t_j}}
 -e^{4\tilde{K}Y^j_{t_{j-1}}}\big]+e^{8\tilde{K}J}(t_j -t_{j-1})\Big\} \\
 &\le &C+C \sum_{j=1}^n
 E\big[e^{4\tilde{K}X_{t_j}}-e^{4\tilde{K}X_{t_{j-1}}}\big]= C+ C
 E\big[e^{4\tilde{K}X_{t_n}}-e^{4\tilde{K}X_{t_0}}\big] \le C,
 \eeas
where we applied (\ref{bmo1}) and (\ref{usis1}) with
$\tilde{K} \dfnn \frac{1}{2}\vee k(J+1) \vee \ell(J)$ to derive the
second inequality and the third inequality is due to the fact that
$Y^j_{t_{j-1}} = \cE^g_{t_{j-1},t_j}[X_{t_j}] \ge X_{t_{j-1}}$. The
proof is now complete. \qed

\ms

With the above upcrossing inequality, we can discuss the continuity
of the quadratic $g$-sub(super)martingales.
 \begin{cor}
 \label{cgsm}
If $X$ is a $g$-submartingale\(resp. $g$-supermartingale\),
then for any denumerable dense subset $\cD$ of $[0,T]$, it holds
$P$-a.s. that \beas \underset{r \nearrow t,\, r \in \cD}{\lim}
X_r\mbox{ exists for any }t \in (0,T]\mbox{ and }\underset{r
\searrow t, \, r \in \cD}{\lim} X_r\mbox{ exists for any }t \in [0,T).
\eeas
\end{cor}

{\it Proof.} If $X$ is a $g$-supermartingale, then $-X$ is
correspondingly a $g^{-}$-submartingale with $g^{-}$ defined in
(\ref{gneg}).
Hence, it suffices to assume that $X$ is a
$g$-submartingale. Let $\{\cD_n\}_{n \in \hN}$ be an increasing
sequence  of finite subsets of $\cD$ such that $\underset{n}
{\cup}\cD_n=\cD$. For any two real numbers $a<b$, Theorem
\ref{ucism} and Jensen's Inequality imply that:
 \beas
 \wt{C}&\dfnn&1+\frac{\|X\|_\infty+k(J+1)T+|a|}{b-a}\\
 &\geq & 1+E\bigg[U^b_a(\wt{X},\cD_n)\exp\Big\{\int_0^{t_n} \beta_\cD(s) dB_s
 - \frac{1}{2} \int_0^{t_n} |\beta_\cD(s)|^2 ds \Big\}\bigg] \\
 &=& E\bigg[\big(1+U^b_a(\wt{X},\cD_n)\big)\exp\Big\{\int_0^{t_n} \beta_\cD(s) dB_s
 - \frac{1}{2} \int_0^{t_n} |\beta_\cD(s)|^2 ds \Big\}\bigg] \\
 &\geq& \exp
 \bigg\{E\Big[\ln{\big(1+U^b_a(\wt{X},\cD_n)\big)}+\int_0^{t_n}
 \beta_\cD(s) dB_s - \frac{1}{2} \int_0^{t_n} |\beta_\cD(s)|^2 ds
 \Big] \bigg\},
 \eeas
from which one can deduce that
 \beas
 E\Big[\ln{\big(1+U^b_a(\wt{X},\cD_n)\big)}\Big] \le \ln\wt{C}+\frac12
 + \|\beta_\cD  \|^2_{L^2_\bF([0,t_n]; \hR^d)} \le C',
 \eeas
where $C'$ is a constant independent of the choice of $\cD$. Since $
U^b_a(\wt{X},\cD_n)\nearrow U^b_a\big(\wt{X},\cD \big)$ as $\cD_n
\nearrow \cD$, Monotone Convergence Theorem implies that
$\ln{\big(1+U^b_a\big(\wt{X},\cD \big)\big)}$ is integrable, thus
$U^b_a\big(\wt{X},\cD \big)< \infty$, $\pas$ Then a classical
argument yields the conclusion for $\wt{X}$, thus for $X$. The proof
is now complete. \qed

\no{\bf Acknowledgment.} We would like to express our sincere gratitude to the anonymous referee for
his/her careful reading of the original manuscript and many valuable suggestions, which helped us to
improve the quality of the paper significantly.

\bibliographystyle{plain}
\bibliography{ma_yao}

\end{document}